\documentclass[a4paper,11pt]{article}
\usepackage{amsmath}
\usepackage{amsfonts}
\usepackage{amssymb}
\usepackage{latexsym}
\usepackage{epsfig}
\usepackage{graphicx}
\usepackage{oldgerm}
\usepackage{theorem}

\def\R{\mathbb{R}}
\def\W{\mathbb{W}}
\def\N{\mathbb{N}}
\def\C{\mathbb{C}}

\def\mbK{\mathbb{K}}
\def\P{\mathbb{P}}

\def\X{\mib{X}}
\def\x{\mib{x}}
\def\y{\mib{y}}
\def\z{\mib{z}}
\def\valpha{\mib{\alpha}}
\def\v{\mib{v}}
\def\w{\mib{w}}
\def\u{\mib{u}}
\def\a{\mib{a}}
\def\b{\mib{b}}
\def\c{\mib{c}}

\def\I{{\bf{I}}}
\def\q{{\bf{q}}}

\def\mM{\mathfrak{M}}

\def\fd{\stackrel{\rm f.d.}{=}}

\def\im{{\rm i}}
\def\supp{{\rm supp}\ }
\def\bE{{\bf E}}
\def\S{{\cal S}}
\def\c{{\rm C}}
\def\mbK{\mathbb{K}}

\theorembodyfont{\itshape}
\newtheorem{thm}{Theorem}[section]
\newtheorem{lem}[thm]{Lemma}
\newtheorem{cor}[thm]{Corollary}

\newcommand{\mib}[1]{\mbox{\boldmath $#1$}}
\newcommand{\SSC}[1]{\section{#1}\setcounter{equation}{0}}
\newcommand{\qed}{\hbox{\rule[-2pt]{3pt}{6pt}}}

\begin{document}

\title{
Characteristic Polynomials of Random Matrices\\
and Noncolliding Diffusion Processes}
\author{
Makoto Katori \\
Department of Physics,
Chuo University \\
{\small Kasuga, Bunkyo-ku, Tokyo 112-8551, Japan} \\
e-mail: katori@phys.chuo-u.ac.jp
}
\date{17 December 2015}
\pagestyle{plain}
\maketitle
\begin{abstract}
We consider the noncolliding Brownian motion (BM)
with $N$ particles starting from the eigenvalue distribution
of Gaussian unitary ensemble (GUE) of $N \times N$
Hermitian random matrices with variance $\sigma^2$.
We prove that this process is equivalent with the
time shift $t \to t+\sigma^2$ of the noncolliding BM
starting from the configuration in which
all $N$ particles are put at the origin.
In order to demonstrate nontriviality of such
equivalence for determinantal processes,
we show that, even from its special consequence,
determinantal expressions are derived for 
the ensemble averages of products of 
characteristic polynomials of random matrices
in GUE.
Another determinantal process, 
noncolliding squared Bessel process with
index $\nu >-1$, is also studied
in parallel with the noncolliding BM
and corresponding results for characteristic polynomials
are given for random matrices in the chiral GUE as well as
in the Gaussian ensembles of class C and class D.

\end{abstract}

\normalsize

\SSC{Introduction}

We consider $N$-particle systems of the one-dimensional
standard Brownian motions (BMs),
$\X(t)=(X_1(t), X_2(t), \dots, X_N(t)), t \geq 0$, 
and of the squared Bessel processes (BESQ) with
index $\nu > -1$,
$\X^{(\nu)}(t)=(X^{(\nu)}_1(t), X^{(\nu)}_2(t), \dots,
X^{(\nu)}_N(t)), t \geq 0$,
both {\it conditioned never to collide with each other},
$N \in \N \equiv \{1,2,3, \dots\}$. 
The former process, which is called the
{\it noncolliding BM} \cite{KT07b},
solves the following set of stochastic differential
equations (SDEs)
\begin{equation}
dX_j(t)=dB_j(t)+ 
\sum_{1 \leq k \leq N, k \not= j}
\frac{dt}{X_j(t)-X_k(t)},
\quad
1 \leq j \leq N, \quad t \geq 0,
\label{eqn:noncollBM}
\end{equation}
with independent one-dimensional standard
BMs $\{B_j(t) \}_{j=1}^{N}$,
and the latter process,
the {\it noncolliding BESQ} \cite{KT11}, does the following set of SDEs
\begin{eqnarray}
dX^{(\nu)}_j(t) &=& 2 \sqrt{X^{(\nu)}_j(t)} 
d \widetilde{B}_j(t) + 2 (\nu+1) dt
\nonumber\\
&& + 4 X^{(\nu)}_j(t)
\sum_{1 \leq k \leq N, k \not= j}
\frac{dt}{X^{(\nu)}_j(t)-X^{(\nu)}_k(t)},
\quad 1 \leq j \leq N, \quad t \geq 0,
\label{eqn:noncollBESQ}
\end{eqnarray}
where $\{\widetilde{B}_j(t)\}_{j=1}^{N}$
are independent one-dimensional standard BMs
different from $\{B_j(t)\}_{j=1}^{N}$
and, if $-1 < \nu < 0$,
the reflection boundary condition is assumed at
the origin.
(See \cite{FNH99,BHL99,KO01,KTNK03,Nag03,TW04,KT07a,TW07,
SMCR08,KMW09,RS10,Ols10,KMW10} 
for related interacting particle systems.) 
Let $\R$ be the collection of all real numbers and
$\R_+=\{x \in \R: x \geq 0\}$, 
and consider subsets of the $N$-dimensional real space
$\R^{N}$,
$\W_N^{\rm A}=\{\x=(x_1, x_2, \dots, x_N) \in \R^N :
x_1 < x_2 < \cdots < x_N\}$,
and 
$\W_N^{+}=\{\x=(x_1, \dots, x_N) \in \R_+^N: 
x_1 < \cdots < x_N\}$.
The former is called the Weyl chambers of
types A$_{N-1}$.
We can prove that, provided $\X(0) \in \W_N^{\rm A}$
and $\X^{(\nu)}(0) \in \W_N^{+}$,
then the SDEs (\ref{eqn:noncollBM})
and (\ref{eqn:noncollBESQ}) guarantee
that with probability one
$\X(t) \in \W_N^{\rm A}$, and
$\X^{(\nu)}(t) \in \W_N^{+}$
for all $t > 0$.
In both processes, at any positive time $t > 0$ 
there is no {\it multiple point}
at which coincidence of particle positions $X_j(t)=X_{k}(t)$ or
$X^{(\nu)}_j(t) = X^{(\nu)}_k(t)$ 
for $j \not= k$ occurs.
It is the reason why these processes are called
noncolliding diffusion processes \cite{KT_Sugaku}.
We can consider them,
however, starting from initial configurations
with multiple points.
An extreme example is the initial configuration
in which all $N$ particles are put at the origin.
In order to describe configurations
with multiple points we represent each
particle configuration by a sum of delta measures
in the form 
$\xi(\cdot)=\sum_{j=1}^N \delta_{x_j}(\cdot)$,
where with given $y \in \R$, $\delta_y(\cdot)$
denotes the delta measure such that
$\delta_{y}(x)=1$ for $x=y$ and
$\delta_y(x)=0$ for $x \not= y$.
Note that, by this definition, for $A \subset \R$, 
$\xi(A)=\int_{A} \xi(dx)=\sum_{1 \leq j \leq N: x_j \in A} 1
=$ the number of particles included in $A$.
(The above mentioned example is then expressed by
$\xi(\cdot)=N \delta_0(\cdot)$,
which means that the origin is the multiple point
with all $N$ particles.)
For a given total number of particles $N \in \N$, 
we write the configuration spaces as
$\mM_N=\{\xi(\cdot)=\sum_{j=1}^N \delta_{x_j}(\cdot) :
x_j \in \R, 1 \leq j \leq N \}$
and
$\mM_N^{+}=\{\xi(\cdot)=\sum_{j=1}^{N} \delta_{x_j}(\cdot) :
x_j \in \R_+, 1 \leq j \leq N \}$.
We consider the noncolliding BM and the 
noncolliding BESQ as $\mM_N$-valued and
$\mM_N^+$-valued processes and write them as
$$
\Xi(t, \cdot)=\sum_{j=1}^N \delta_{X_j(t)}(\cdot),
\quad
\Xi^{(\nu)}(t, \cdot)=\sum_{j=1}^N
\delta_{X^{(\nu)}_j(t)}(\cdot),
\quad t \geq 0,
$$
respectively. 
The probability law of $\Xi(t, \cdot)$
starting from a fixed configuration $\xi \in \mM_N$
is denoted by $\P^{\xi}$ and that of $\Xi^{(\nu)}(t, \cdot)$
from $\xi \in \mM_N^{+}$ by $\P_{\nu}^{\xi}$,
and the noncolliding diffusion processes
specified by initial configurations 
are expressed by
$(\Xi(t), t \in [0, \infty), \P^{\xi})$ and
$(\Xi^{(\nu)}(t), t \in [0, \infty), \P^{\xi}_{\nu}), \nu > -1$,
respectively.
We set 
$
\mM_{N,0}= \{ \xi \in \mM_N : \xi(\{x\})\le 1 \mbox { for any }  x\in\R \}$,
and
$\mM^+_{N,0}= \{ \xi \in \mM_N^{+} : \xi(\{x\})\le 1 \mbox { for any }  
x \in \R_+ \}$,
which denote collections of configurations
without any multiple points.

In order to dynamically simulate the random matrix 
ensemble called the Gaussian unitary ensemble (GUE),
Dyson considered the $N \times N$ Hermitian matrix-valued
BM and showed that its eigenvalue process satisfies
the SDEs given by (\ref{eqn:noncollBM}) \cite{Dys62}.
This eigenvalue process is called 
Dyson's BM model with parameter $\beta=2$
or simply Dyson's model \cite{Spo87,KT10}.
The equivalence between Dyson's model
and the noncolliding BM, $(\Xi(t), t \in [0, \infty), \P^{\xi})$, 
implies that the random matrix theory \cite{Meh04,For10}
is useful to classify and analyze noncolliding diffusion processes
\cite{KT04,KT_Sugaku}. 
In particular, if the initial configuration is given by
$\xi=N \delta_0$, this equivalence concludes that, 
for any $t > 0$,
the particle distribution of $\Xi(t)$ is equal to the
eigenvalue distribution of $N \times N$
random matrices in the GUE with variance $t$.
(Note that in the usual GUE the mean is set to be zero.)
Here the probability density function (pdf) of
the GUE eigenvalues with variance $\sigma^2$
is given by
\begin{equation}
\mu_{N, \sigma^2}(\xi)
=\frac{\sigma^{-N^2}}{C_N}
\exp \left( - \frac{|\x|^2}{2 \sigma^2} \right)
h_N(\x)^2,
\label{eqn:GUE1}
\end{equation}
$\xi=\sum_{j=1}^{N} \delta_{x_j} \in \mM_N, 
x_1 \leq x_2 \leq \cdots \leq x_N$,
where $C_N=(2 \pi)^{N/2} \prod_{j=1}^{N} \Gamma(j)$
with the gamma function 
$\Gamma(z)=\int_0^{\infty} e^{-u} u^{z-1} du$,
$|\x|^2=\sum_{j=1}^{N} x_j^2$,
and 
\begin{equation}
h_{N}(\x)=\prod_{1 \leq j< k \leq N}(x_k-x_j).
\label{eqn:Vandermonde}
\end{equation}
The expectation of a measurable function $F$ of
a random variable $\Xi \in \mM_N$ with respect to
(\ref{eqn:GUE1}) is given by
\begin{equation}
\bE_{N, \sigma^2} [F(\Xi)]
=\int_{\W_N^{\rm A}} F(\xi) \mu_{N, \sigma^2}(\xi) d \x
=\frac{1}{N!} \int_{\R^N} F(\xi) \mu_{N, \sigma^2}(\xi) d \x
\label{eqn:ENs}
\end{equation}
with setting $\xi=\sum_{j=1}^{N} \delta_{x_j},
\x=(x_1, \dots, x_N)$,
where $d \x=\prod_{j=1}^N dx_j$.

Assume $\nu \in \N_0 \equiv \N \cup \{0\}$
and let $M(t)$ be the $(N+\nu) \times N$ 
complex matrix-valued BM.
Then the $N \times N$ matrix-valued process
$L(t)=M(t)^{\dagger} M(t), t \geq 0$,
where $M(t)^{\dagger}$ denotes the Hermitian conjugate
of $M(t)$, is called the Laguerre process
or the complex Wishart process \cite{Bru91}.
The matrix $L(t)$ is Hermitian and positive definite,
and K\"onig and O'Connell proved that
the eigenvalue process of $L(t)$ satisfies the SDEs
given by (\ref{eqn:noncollBESQ}) \cite{KO01}.
Again by the random matrix theory \cite{Meh04,For10},
this equivalence concludes that the particle distribution
of $(\Xi^{(\nu)}(t), t \in [0, \infty), \P_{\nu}^{N \delta_0})$
at any given time $t >0$
is equal to the distribution of squares of singular values
of $(N+\nu) \times N$ complex random matrices in the
chiral Gaussian unitary ensemble (chGUE)
with index $\nu \in \N_0$ and variance $t$.
Here the pdf of squares of singular values in
chGUE with index $\nu \in \N_0$ and variance $\sigma^2$ is given by
\cite{VZ93,Ver94}
\begin{equation}
\mu_{N, \sigma^2}^{(\nu)}(\xi)
=\frac{\sigma^{-2N(N+\nu)}}{C_N^{(\nu)}}
\prod_{j=1}^{N} ( x_j^{\nu} e^{-x_j/2 \sigma^2} )
h_N(\x)^2,
\label{eqn:chGUE1}
\end{equation}
$\xi=\sum_{j=1}^{N} \delta_{x_j} \in \mM_N^+,
0 \leq x_1 \leq \cdots \leq x_N$,
where $C_N^{(\nu)}=2^{N(N+\nu)} \prod_{j=1}^{N}
\Gamma(j) \Gamma(j+\nu)$.
The expectation of a measurable function $F$ of
a random variable $\Xi \in \mM^+_N$ with respect to
(\ref{eqn:chGUE1}) is given by
\begin{equation}
\bE^{(\nu)}_{N, \sigma^2} [F(\Xi)]
=\int_{\W_N^{+}} F(\xi) \mu_{N, \sigma^2}^{(\nu)}(\xi) d \x
=\frac{1}{N!} \int_{\R_+^N} F(\xi) \mu_{N, \sigma^2}^{(\nu)}(\xi) d \x
\label{eqn:ENsn}
\end{equation}
with setting $\xi=\sum_{j=1}^{N} \delta_{x_j},
\x=(x_1, \dots, x_N)$.

Let ${\cal H}(2N)$ be the space of $2N \times 2N$ Hermitian
matrices and $\mathfrak{sp}(2N, \C)$ and $\mathfrak{so}(2N, \C)$
be the spaces of $2N \times 2N$ complex matrices
representing the symplectic Lie algebra
and the orthogonal Lie algebra, respectively.
Altland and Zirnbauer introduced
the Gaussian random matrix ensembles
for the elements in 
${\cal H}^{\rm C}(2N)=\mathfrak{sp} \cap {\cal H}(2N)$
and in ${\cal H}^{\rm D}(2N)=\mathfrak{so} \cap {\cal H}(2N)$,
which are called the Gaussian ensembles of class C and class D,
respectively. 
The eigenvalues of matrices both in the class C and class D ensembles
are given by $N$ pairs of positive and negative ones
with the same absolute value.
The pdfs of the squares of $N$ positive eigenvalues
are given by (\ref{eqn:chGUE1}) with
$\nu=1/2$ for the class C ensemble
and with $\nu=-1/2$ for the class D ensemble,
when the variances are $\sigma^2$ \cite{AZ96,AZ97}.
If we consider the ${\cal H}^{\rm C}(2N)$-valued BM
and ${\cal H}^{\rm D}(2N)$-valued BM,
the squares of each $N$ positive eigenvalues
satisfy the SDEs (\ref{eqn:noncollBESQ}) with
$\nu=1/2$ and $\nu=-1/2$, respectively \cite{KT04}.
See also \cite{BHL99,KTNK03,Nag03,BFPSW09,KS10}.

In the present paper, we consider the noncolliding BM
whose initial configuration is distributed according
to the pdf (\ref{eqn:GUE1}), denoted by
$(\Xi(t), t \in [0, \infty), \P^{\mu_{N, \sigma^2}})$,
and the noncolliding BESQ starting from the
distribution (\ref{eqn:chGUE1})
with not only $\nu \in \N_0$ but with 
$\nu > -1$ generally, denoted by
$(\Xi^{(\nu)}(t), t \in [0, \infty), \P^{\mu_{N, \sigma^2}^{(\nu)}})$.
We prove that for any $N \in \N, \sigma^2 > 0$ 
the following equalities are established,
\begin{eqnarray}
(\Xi(t), t \in [0, \infty), \P^{\mu_{N, \sigma^2}}) 
&\fd& (\Xi(t+\sigma^2), t \in [0, \infty), \P^{N \delta_0}),
\nonumber\\
(\Xi^{(\nu)}(t), t \in [0, \infty), \P_{\nu}^{\mu^{(\nu)}_{N, \sigma^2}})
&\fd& (\Xi^{(\nu)}(t+\sigma^2), t \in [0, \infty), \P_{\nu}^{N \delta_0}),
\quad \nu > -1,
\label{eqn:equal1}
\end{eqnarray}
where $\fd$ denotes the equivalence in finite dimensional
distributions (see Theorem \ref{Theorem:distributions} and Remark 1).
We would like to emphasize that these equalities are
highly nontrivial and in order to demonstrate it
we show in this paper that,
even from very special consequence of (\ref{eqn:equal1}),
determinantal expressions of ensemble averages
of products of characteristic polynomials
of random matrices are derived.
See \cite{BH00,FW01,MN01,BH01,FS02,FS03,AV03,BS06}
for extensive study
of characteristic polynomials of random matrices,
especially in the connection with the
Riemann zeta function \cite{KS00a,KS00b,HKO00}.

We write the expectations of
measurable function $F$ of $N \times N$ random matrices $\{H\}$ 
in the GUE, 
of $(N+\nu) \times N$ random matrices $\{M\}$ 
in the chGUE 
with $\nu \in \N_0$,
of $2N \times 2N$ random matrices $\{H^{\rm C}\}$
in the class C,
and of $2N \times 2N$ random matrices $\{H^{\rm D}\}$
in the class D as 
$\langle F(H) \rangle_{{\rm GUE}(N, \sigma^2)}$,
$\langle F(M) \rangle_{{\rm chGUE}(N, \nu, \sigma^2)}$, 
$\langle F(H^{\rm C}) \rangle_{{\rm classC}(2N, \sigma^2)}$,
and
$\langle F(H^{\rm D}) \rangle_{{\rm classD}(2N, \sigma^2)}$,
respectively, where
$\sigma^2$ denote the variances of these four kinds of
Gaussian ensembles.
Then for $m \in \N$, $\valpha \in \C^{m}$
the ensemble averages of $m$-product of 
characteristic polynomials of random matrices
are defined as
\begin{equation}
M_{\rm GUE}(m, \valpha; N, \sigma^2)
\equiv \left\langle
\prod_{n=1}^{m} \det (\alpha_n I_N - H) 
\right\rangle_{{\rm GUE}(N, \sigma^2)}
= \left\langle
\prod_{n=1}^{m} \prod_{j=1}^{N} (\alpha_n -\lambda_j) 
\right\rangle_{{\rm GUE}(N, \sigma^2)}, 
\label{eqn:PGUE1}
\end{equation}
\begin{eqnarray}
M^{(\nu)}_{\rm chGUE}(m, \valpha; N, \sigma^2)
&\equiv& \left\langle
\prod_{n=1}^{m} \det (\alpha_n I_N - M^{\dagger} M) 
\right\rangle_{{\rm chGUE}(N,\nu,\sigma^2)}
\nonumber\\
&=& \left\langle
\prod_{n=1}^{m} \prod_{j=1}^{N} (\alpha_n -\kappa_j^2) 
\right\rangle_{{\rm chGUE}(N,\nu,\sigma^2)}, 
\quad \nu \in \N_0,
\label{eqn:PchGUE1}
\end{eqnarray} 
and for $\sharp=\mbox{C and D}$
\begin{eqnarray}
&& M_{{\rm class}\sharp}(m, \valpha; 2N, \sigma^2)
\equiv
\left\langle \prod_{n=1}^m \det (\alpha_n I_{2N}-H^{\sharp})
\right\rangle_{{\rm class}\sharp(2N,\sigma^2)}
\nonumber\\
\label{eqn:Pclass}
&& 
= \left\langle \prod_{n=1}^m \prod_{j=1}^{N}
(\alpha_n-\varepsilon_j)(\alpha_n+\varepsilon_j)
\right\rangle_{{\rm class}\sharp(2N,\sigma^2)}
= \left\langle \prod_{n=1}^m \prod_{j=1}^{N}
(\alpha_n^2-\varepsilon_j^2)
\right\rangle_{{\rm class}\sharp(2N,\sigma^2)},
\end{eqnarray}
where $I_{\ell}$ denotes the $\ell \times \ell$ unit matrix,
$(\lambda_1, \dots, \lambda_N)$ are the eigenvalues of $H$,
$(\kappa_1^2, \dots, \kappa_N^2)$ are the
eigenvalues of $M^{\dagger} M$,
$(\varepsilon_1, \dots, \varepsilon_N, -\varepsilon_1,
\dots, -\varepsilon_N)$ are
the eigenvalues forming 
``particle-hole pairing" of $H^{\sharp}, \sharp=$ C or D.
For (\ref{eqn:PchGUE1}) 
remark that each $M$ in chGUE has such a 
singular value decomposition that
$M=U^{\dagger} K V$, where
$U \in {\rm U}(N+\nu), V \in {\rm U}(N)$, 
$$
  K = \left( \begin{array}{c}
  \widehat{K} \cr O 
  \end{array} \right)
\quad
\mbox{with} \quad
\widehat{K}={\rm diag}
\{\kappa_{1}, \kappa_{2}, \cdots, \kappa_{N}\}, \quad
(\kappa_1, \dots, \kappa_N) \in \W_N^{+}
$$
and $O$ is the $\nu \times N$ zero matrix.
The diagonal elements $(\kappa_1, \dots, \kappa_N)$
of $\widehat{K}$
are called the singular values of 
rectangular matrix $M$. 
Since $M^{\dagger} M=V^{\dagger} K^{\rm T} K V$,
the eigenvalues of $M^{\dagger} M$ are
squares of singular values
$(\kappa_1^2, \dots, \kappa_N^2)$ \cite{Hua63}.

We use the convention such that
$$
\prod_{x\in\xi}f(x) =\exp
\left\{\int_\R \xi(dx) \log f(x) \right\}
=\prod_{x \in \supp \xi}f(x)^{\xi(\{x\})}
$$
for $\xi\in \mM_N$ and a function $f$ on $\R$,
where $\supp \xi = \{x \in \R : \xi(\{x\}) > 0\}$.
Then (\ref{eqn:PGUE1}) is given by
\begin{equation}
M_{\rm GUE}(m, \valpha; N, \sigma^2)
= \bE_{N, \sigma^2} \left[
\prod_{n=1}^{m} \prod_{X \in \Xi}(\alpha_n-X) \right]
\label{eqn:M1}
\end{equation}
with (\ref{eqn:ENs}).
And if we define
\begin{equation}
M^{(\nu)}(m, \valpha; N, \sigma^2)
= \bE^{(\nu)}_{N, \sigma^2} \left[
\prod_{n=1}^{m} 
\prod_{X \in \Xi}(\alpha_n-X) \right]
\label{eqn:PchGUE2}
\end{equation}
with (\ref{eqn:ENsn}) for $\nu > -1$, 
(\ref{eqn:PchGUE1}) and (\ref{eqn:Pclass}) 
are given as
\begin{eqnarray}
\label{eqn:M2}
&& M_{\rm chGUE}^{(\nu)}(m, \valpha; N, \sigma^2)
=M^{(\nu)}(m, \valpha; N, \sigma^2),
\quad \nu \in \N_0, \\
\label{eqn:M3}
&& M_{\rm classC}(m, \valpha; 2N, \sigma^2)
=M^{(1/2)}(m, \valpha^{\langle 2 \rangle}; N, \sigma^2), \\
\label{eqn:M4}
&& M_{\rm classD}(m, \valpha; 2N, \sigma^2)
=M^{(-1/2)}(m, \valpha^{\langle 2 \rangle}; N, \sigma^2),
\end{eqnarray}
where for $\valpha=(\alpha_1, \dots, \alpha_m) \in \C^m$,
$\valpha^{\langle 2 \rangle} \equiv
(\alpha_1^2, \dots, \alpha_m^2) \in \C^m$.

We will show that, from the equalities (\ref{eqn:equal1}),
two sets of determinantal expressions are derived for (\ref{eqn:M1}),
(\ref{eqn:M2})-(\ref{eqn:M4}).

\SSC{Preliminaries and Main Results}

We define
\begin{equation}
p(t, y|x)
=  \left\{ \begin{array}{ll}
\displaystyle{
\frac{1}{\sqrt{2 \pi |t|}} 
\exp \left\{ - \frac{(x-y)^2}{2t} \right\}},
& \quad t \in \R \setminus \{0\},  \cr
\delta(y-x),
& \quad t=0, 
\end{array} \right.
\label{eqn:p}
\end{equation}
for $x, y \in \C$.
For $z \in \C, \nu > -1$, we define $z^{\nu}$ 
to be $\exp(\nu \log z)$, where the argument of $z$
is given its principal value;
$$
z^{\nu}=\exp \Big[ \nu \Big\{
\log |z| + \sqrt{-1} {\rm arg} (z) \Big\} \Big], \quad
-\pi < {\rm arg} (z) \leq \pi.
$$
For $\nu > -1, y \in \C$, we set
\begin{equation}
p^{(\nu)}(t, y|x)
=  \left\{ \begin{array}{ll}
\displaystyle{
\frac{1}{2|t|} \left( \frac{y}{x} \right)^{\nu/2}
\exp \left( - \frac{x+y}{2t} \right)
I_{\nu} \left( \frac{\sqrt{xy}}{t} \right)},
& \quad t \in \R \setminus \{0\}, x \in \C \setminus \{0\}, \cr
\displaystyle{
\frac{y^{\nu}}{(2 |t|)^{\nu+1} \Gamma(\nu+1)} e^{-y/2t}},
& \quad t \in \R \setminus \{0\}, x=0, \cr
& \cr
\delta(y-x),
& \quad t=0, x \in \C,
\end{array} \right.
\label{eqn:pnu}
\end{equation}
where $I_{\nu}(z)$ is the modified Bessel function
of the first kind defined by \cite{Wat44,AAR99}
$$
I_{\nu}(z) = \sum_{n=0}^{\infty} 
\frac{1}{\Gamma(n+1) \Gamma(n+1+\nu)}
\left( \frac{z}{2} \right)^{2n+\nu}.
$$

The following equality holds,
\begin{eqnarray}
\label{eqn:CK1}
&& \int_{\R} dy \, p(s,z|y) p(t,y|x) = p(s+t, z|x), \\
\label{eqn:CK2}
&& \int_{\R+} dy \, p^{(\nu)}(s,z|y) 
p^{(\nu)}(t,y|x) = p^{(\nu)}(s+t, z|x)
\end{eqnarray}
for $s, t \geq 0, x, z \in \C$.
The former is confirmed just performing the Gaussian integral
and the latter is proved by using Weber's second exponential
integral of the Bessel functions \cite{Wat44}
with appropriate analytic continuation \cite{KT11}.
When $t \geq 0, x, y \in \R$, $p(t, y|x)$
gives the transition probability density of the one-dimensional
standard BM from $x$ to $y$ during time period $t$,
and when $t \geq 0, x, y \in \R_+$,
$p^{(\nu)}(t, y|x)$ gives that of the BESQ with index $\nu>-1$
from $x$ to $y$ during time period $t$
(if $-1 < \nu < 0$, a reflection wall is put at the origin).
The equalities (\ref{eqn:CK1}) with $x, z \in \R$
and (\ref{eqn:CK2}) with $x, z \in \R_+$ are the
Chapman-Kolmogorov equations
of these one-dimensional diffusion processes \cite{RY98}.
The extended versions of Chapmann-Kolmogorov equations,
{\it i.e.}, (\ref{eqn:CK1}) with $x \in \C \setminus \R$
and (\ref{eqn:CK2}) with $x \in \R \setminus \R_+$ will have
probability-theoretical interpretations related with martingales \cite{Kat14}.

We introduce the Karlin-McGregor determinants \cite{KM59}
\begin{eqnarray}
f(t, \y|\x) &=& \det_{1 \leq j, k \leq N}
[p(t, y_j|x_k)], \qquad \x, \y \in \W_N^{\rm A},
\nonumber\\
f^{(\nu)}(t, \y|\x) &=& \det_{1 \leq j, k \leq N}
[p^{(\nu)}(t, y_j|x_k)], \qquad \x, \y \in \W_N^{+},
\quad \nu > -1,
\label{eqn:f1}
\end{eqnarray}
$t \geq 0$.

For given $\xi \in \mM_N, N \in \N$,
we write $\xi(\cdot)=\sum_{j=1}^N \delta_{a_j}(\cdot)$
with $a_1 \leq a_2 \leq \cdots \leq a_N$.
Then we set
$$
\xi_0(\cdot) \equiv 0 \quad
\mbox{and} \quad
\xi_n(\cdot)=\sum_{j=1}^{n} \delta_{a_j}(\cdot),
\quad 1 \leq n \leq N.
$$
For $t > 0, x \in \C$, we define
\cite{BK05,KT10}
\begin{eqnarray}
\label{eqn:phi1a}
\phi_n^{(+)}(t, x; \xi)
&=& \frac{1}{2 \pi \im} 
\oint_{\c(\xi_{n+1})} ds \,
p(t, x|s) \frac{1}{\prod_{a \in \xi_{n+1}}(s-a)},
\\
\label{eqn:phi1b}
\phi_n^{(-)}(t, x; \xi)
&=& \int_{\R} ds \, 
p(-t, \im s|x) \prod_{a \in \xi_n} (\im s-a),
\quad \xi \in \mM_N,
\end{eqnarray}
and \cite{DF08,KT11}
\begin{eqnarray}
\label{eqn:phinu1a}
\phi_n^{(\nu,+)}(t, x; \xi)
&=& \frac{1}{2 \pi \im} 
\oint_{\c(\xi_{n+1})} ds \,
p^{(\nu)}(t, x|s) \frac{1}{\prod_{a \in \xi_{n+1}}(s-a)},
\\
\label{eqn:phinu1b}
\phi_n^{(\nu,-)}(t, x; \xi)
&=& \int_{\R_-} ds \, 
p^{(\nu)}(-t, s|x) \prod_{a \in \xi_n} (s-a),
\quad \xi \in \mM^+_N,
\end{eqnarray}
$n=0,1, \dots, N-1$,
where $\im =\sqrt{-1}$,
for $\zeta \in \mM_{\ell}, \ell \in \N$, 
$\c(\zeta)$ denotes a closed contour
on the complex plane $\C$ encircling the points in $\supp \zeta$
on the real line $\R$ once in the positive direction, 
and $\R_-=\{x \in \R: x \leq 0\}$. 
And for $t \geq 0$, let
\begin{eqnarray}
&& h^{(\pm)}_N(t, \y; \xi)
= \det_{1 \leq j, k \leq N} 
[ \phi^{(\pm)}_{j-1}(t, \y_k; \xi)],
\quad \y \in \W_N^{\rm A}, \quad \xi \in \mM_N,
\nonumber\\
&& h^{(\nu, \pm)}_N(t, \y; \xi)
= \det_{1 \leq j, k \leq N} 
[ \phi^{(\nu, \pm)}_{j-1}(t, \y_k; \xi)],
\quad \y \in \W_N^{+}, \quad \xi \in \mM_N^+.
\label{eqn:hpm}
\end{eqnarray}
Since $\phi_n^{(-)}(t, x; \xi)$ and
$\phi_n^{(\nu,-)}(t, x; \xi)$ are monic polynomials
of $x$ of degree $n$, 
\begin{equation}
h^{(-)}_N(t, \y; \xi)
=h^{(\nu, -)}_N(t, \y; \xi)
=h_N(\y),
\label{eqn:h-}
\end{equation}
which are independent of $t$ and $\xi$.
On the other hand, the following equalities are
proved (Lemma 3.1 in \cite{KT10}
and Lemma 3.4 in \cite{KT11}).
For any $t \geq 0, \y \in \W_N^{\rm A},
\xi=\sum_{j=1}^{N} \delta_{x_j} \in \mM_N$
with $x_1 \leq x_2 \leq \cdots \leq x_N$,
\begin{equation}
\frac{f(t, \y|\x)}{h_N(\x)}
=h^{(+)}_N(t, \y; \xi),
\label{eqn:fhh1}
\end{equation}
and for any $t \geq 0, \y \in \W_N^{+},
\xi=\sum_{j=1}^{N} \delta_{x_j} \in \mM_N^+$
with $0 \leq x_1 \leq x_2 \leq \cdots \leq x_N$,
\begin{equation}
\frac{f^{(\nu)}(t, \y|\x)}{h_N(\x)}
=h^{(\nu,+)}_N(t, \y; \xi),
\label{eqn:fhh2}
\end{equation}
where, if some of the $x_j$'s coincide, the LHS 
of (\ref{eqn:fhh1}) and (\ref{eqn:fhh2}) are
interpreted by using l'H\^opital's rule.

For any $M \in \N$ and any increasing time-sequence
$0 < t_1 < \cdots < t_M < \infty$,
the multitime probability density of 
$(\Xi(t), t \in [0, \infty), \P^{\xi})$
is given by
\begin{eqnarray}
&& p^{\xi}(t_1, \xi^{(1)}; \dots; t_M, \xi^{(M)})
\nonumber\\
&& \qquad = h^{(-)}_{N}(t_M, \x^{(M)}; \xi)
\prod_{m=1}^{M-1} f(t_{m+1}-t_m, \x^{(m+1)}|\x^{(m)})
h^{(+)}_{N}(t_1, \x^{(1)}, \xi)
\label{eqn:pxi}
\end{eqnarray}
with $\xi=\sum_{j=1}^{N} \delta_{a_j} \in \mM_N$,
$a_1 \leq a_2 \leq \dots \leq a_N$
for the initial configuration and
$\xi^{(m)}=\sum_{j=1}^N \delta_{x^{(m)}_j} \in \mM_{N,0}$,
$\x^{(m)}=(x^{(m)}_1, \dots, x^{(m)}_N) \in \W_N^{\rm A}$
for configurations at times $t_m, 1 \leq m \leq M$ \cite{KT10}.
Similarly, 
the multitime probability density of 
$(\Xi^{(\nu)}(t), t \in [0, \infty), \P_{\nu}^{\xi})$
is given by
\begin{eqnarray}
&& p_{\nu}^{\xi}(t_1, \xi^{(1)}; \dots; t_M, \xi^{(M)})
\nonumber\\
&& \qquad  = h^{(\nu,-)}_{N}(t_M, \x^{(M)}; \xi)
\prod_{m=1}^{M-1} f^{(\nu)}(t_{m+1}-t_m, \x^{(m+1)}|\x^{(m)})
h^{(\nu,+)}_{N}(t_1, \x^{(1)}; \xi)
\label{eqn:pxinu}
\end{eqnarray}
with $\xi=\sum_{j=1}^{N} \delta_{a_j} \in \mM_N^+$,
$0 \leq a_1 \leq a_2 \leq \dots \leq a_N$
for the initial configuration and
$\xi^{(m)}=\sum_{j=1}^N \delta_{x^{(m)}_j} \in \mM^+_{N,0}$,
$\x^{(m)}=(x^{(m)}_1, \dots, x^{(m)}_N) \in \W_N^{+}$
for configurations at times $t_m, 1 \leq m \leq M$ \cite{KT11}.
By definitions (\ref{eqn:phi1a}) and 
(\ref{eqn:phi1b}), we can see that
\begin{eqnarray}
\phi^{(+)}_n(t, x; N \delta_0)
&=& t^{-(n+1)/2} \frac{2^{-n/2}}{n! \sqrt{2 \pi}}
H_n \left( \frac{x}{\sqrt{2t}} \right) e^{-x^2/2t},
\nonumber\\
\phi^{(-)}_n(t, x; N \delta_0)
&=& t^{n/2} 2^{-n/2} H_n \left( \frac{x}{\sqrt{2t}} \right),
\quad 0 \leq n \leq N-1,
\label{eqn:Hermite1}
\end{eqnarray}
where $H_n(x)$ is the Hermite polynomial of degree $n$,
\begin{eqnarray}
H_n(x) &=& n ! \sum_{k=0}^{[n/2]}
(-1)^k \frac{(2x)^{n-2k}}{k ! (n-2k)!}
\nonumber\\
&=& 2^{n/2}
\int_{\R} dy \,
\frac{e^{-y^2/2}}{\sqrt{2 \pi}} (\im y + \sqrt{2} x)^n
= \frac{n!}{2 \pi \im} \oint_{\c(\delta_0)}
d z \,
\frac{e^{2 z x-z^2}}{z^{n+1}},
\nonumber\\
\label{eqn:Hermite2}
\end{eqnarray}
where $[r]$ denotes the largest integer
that is not greater than $r \in R_+$,
and by definitions (\ref{eqn:phinu1a}) and 
(\ref{eqn:phinu1b}),
\begin{eqnarray}
\phi^{(\nu,+)}_n(t, x; N \delta_0)
&=& t^{-(n+1)} (-1)^n
\frac{2^{-(n+\nu+1)}}{\Gamma(n+\nu+1)}
\left(\frac{x}{t}\right)^{\nu} e^{-x/2t}
L_{n}^{\nu} \left( \frac{x}{2t} \right),
\nonumber\\
\phi^{(\nu,-)}_n(t, x; N \delta_0)
&=& t^n (-1)^n 2^n n ! 
L_{n}^{\nu} \left( \frac{x}{2t} \right),
\label{eqn:Laguerre1}
\end{eqnarray}
where $L^{(\nu)}_n(x)$ is the Laguerre polynomial of degree $n$
with index $\nu$,
\begin{eqnarray}
L_n^{\nu}(x) &=& \sum_{k=0}^{n} (-1)^k
\frac{\Gamma(n+\nu+1) x^k}
{\Gamma(k+\nu+1) (n-k)! k!}
\nonumber\\
&=& \frac{1}{2 \pi \im}
\oint_{\c(\delta_0)} 
\frac{e^{-x u}(1+u)^{n+\nu}}{u^{n+1}} du.
\label{eqn:Laguerre2}
\end{eqnarray}
Then we can prove the following.
\begin{thm}
\label{Theorem:distributions}
For any $N, M \in \N$, any increasing time-sequence
$0 < t_1 < \cdots < t_M < \infty$, and 
any $\sigma^2 >0$,
\begin{equation}
\bE_{N, \sigma^2}
[p^{\Xi}(t_1, \xi^{(1)}; \dots; t_M, \xi^{(M)})]
=p^{N \delta_0}(t_1+\sigma^2, \xi^{(1)};
\dots; t_M+\sigma^2, \xi^{(M)})
\label{eqn:Th1a1}
\end{equation}
with $\xi^{(m)} \in \mM_{N,0}, 1 \leq m \leq M$, and
\begin{equation}
\bE^{(\nu)}_{N, \sigma^2}
[p_{\nu}^{\Xi}(t_1, \xi^{(1)}; \dots; t_M, \xi^{(M)})]
=p_{\nu}^{N \delta_0}(t_1+\sigma^2, \xi^{(1)};
\dots; t_M+\sigma^2, \xi^{(M)})
\label{eqn:Th1a2}
\end{equation}
with $\xi^{(m)} \in \mM^+_{N,0}, 1 \leq m \leq M$.
\end{thm}
\vskip 0.5cm
\noindent{\bf Remark 1.}
When $M$-time probability density of a process
is given for any $M \in \N$ and any time sequence
$0 < t_1 < \cdots < t_M$,
it is said that the finite dimensional distributions
of the process is determined \cite{RY98}.
Eq.(\ref{eqn:Th1a1}) (resp. Eq.(\ref{eqn:Th1a2}))
means that the processes
$(\Xi(t), t \in [0, \infty), \P^{\mu_{N, \sigma^2}})$
and 
$(\Xi(t+\sigma^2), t \in [0, \infty), \P^{N \delta_0})$
(resp. $(\Xi^{(\nu)}(t), t \in [0, \infty), \P_{\nu}^{\mu_{N, \sigma^2}})$
and 
$(\Xi^{(\nu)}(t+\sigma^2), t \in [0, \infty), \P_{\nu}^{N \delta_0})$)
are equivalent in finite dimensional distributions
for any $\sigma^2 >0$,
which is the fact expressed by (\ref{eqn:equal1}).
\vskip 0.5cm

For $\x^{(m)}=(x^{(m)}_1, \dots, x^{(m)}_N) \in \W_N^{\rm A}$ 
or $\x^{(m)} \in \W_N^{+}$ with
$\xi^{(m)}(\cdot)=\sum_{j=1}^{N} \delta_{x^{(m)}_j} (\cdot)$
and $N' \in \{ 1,2, \dots, N\}$, we put
$\x^{(m)}_{N'}=(x^{(m)}_1, \dots, x^{(m)}_{N'})
\in \W_{N'}^{\rm A}$ or $\x^{(m)}_{N'} \in \W_{N'}^{+}$,
$1 \leq m \leq M$.
For a sequence $(N_m)_{m=1}^{M}$ of positive integers 
less than or equal to $N$,
we define the 
$(N_1, \dots, N_{M})$-{\it multitime correlation functions} 
of $(\Xi(t), t \in [0, \infty), \P^{\xi})$ 
and $(\Xi^{(\nu)}(t), t \in [0, \infty), \P^{\xi}_{\nu})$ by
\begin{eqnarray}
&& \rho^{\xi} (t_{1}, \x^{(1)}_{N_1} ; 
\dots; t_M, \x^{(M)}_{N_M}) 
\nonumber\\
&&=
\int_{\prod_{m=1}^{M} \R^{N-N_{m}}}
\prod_{m=1}^{M} \prod_{j=N_{m}+1}^{N} dx_{j}^{(m)}
p^{\xi} (t_1, \xi^{(1)}; \dots; t_M, \xi^{(M)})
\prod_{m=1}^{M}
\frac{1}{(N-N_{m})!},
\nonumber\\
&& \rho^{\xi}_{\nu} (t_{1}, \x^{(1)}_{N_1} ; 
\dots; t_M, \x^{(M)}_{N_M}) 
\nonumber\\
&&=
\int_{\prod_{m=1}^{M} \R_+^{N-N_{m}}}
\prod_{m=1}^{M} \prod_{j=N_{m}+1}^{N} dx_{j}^{(m)}
p_{\nu}^{\xi} (t_1, \xi^{(1)}; \dots; t_M, \xi^{(M)} )
\prod_{m=1}^{M}
\frac{1}{(N-N_{m})!},
\label{eqn:corr}
\end{eqnarray}
respectively.
In the previous papers we have shown that 
for any fixed initial configuration $\xi$ 
the noncolliding BM and the noncolliding BESQ
with finite numbers of particles 
are {\it determinantal processes} in the sense that
any multitime correlation function is given by a determinant
\cite{KT10,KT11}
\begin{eqnarray}
&& \rho^{\xi}
 (t_1,\x^{(1)}_{N_1}; \dots;t_M,\x^{(M)}_{N_M} ) 
=\det_{
\substack{1 \leq j \leq N_{m}, 1 \leq k \leq N_{n} \\ 1 \leq m, n \leq M}
}
[
\mbK^{\xi}(t_m, x_{j}^{(m)}; t_n, x_{k}^{(n)} )]
\nonumber\\
&& \rho^{\xi}_{\nu}
 (t_1,\x^{(1)}_{N_1}; \dots;t_M,\x^{(M)}_{N_M}) 
=\det_{
\substack{1 \leq j \leq N_{m}, 1 \leq k \leq N_{n} \\ 1 \leq m, n \leq M}
}
[
\mbK^{\xi}_{\nu}(t_m, x_{j}^{(m)}; t_n, x_{k}^{(n)} )],
\quad \nu > -1.
\label{eqn:rho1}
\end{eqnarray}
Here the correlation kernels are given by
\begin{eqnarray}
\mbK^{\xi}(s, x; t, y)
&=& \sum_{n=0}^{N-1} \phi_n^{(+)}(s, x; \xi)
\phi_n^{(-)}(t,y; \xi)
- {\bf 1}(s>t) p(s-t, x|y)
\nonumber\\
&=& \int_{\R} \xi(dx') \int_{\R} dy' p(s,x|x') 
\Phi(\xi; x', \im y') p(-t, \im y'|y)
\nonumber\\
\label{eqn:Kxi}
&& - {\bf 1}(s>t) p(s-t, x|y),
\quad \xi \in \mM_{N,0}, \\
\mbK^{\xi}_{\nu}(s, x; t, y)
&=& \sum_{n=0}^{N-1} \phi_n^{(\nu,+)}(s, x; \xi)
\phi_n^{(\nu,-)}(t,y; \xi)
- {\bf 1}(s>t) p^{(\nu)}(s-t, x|y)
\nonumber\\
&=& \int_{\R_+} \xi(dx') \int_{\R_-} dy' p^{(\nu)}(s,x|x') 
\Phi(\xi; x', y') p^{(\nu)}(-t, y'|y)
\nonumber\\
\label{eqn:Kxinu}
&& - {\bf 1}(s>t) p^{(\nu)}(s-t, x|y),
\quad \xi \in \mM_{N,0}^+,
\quad \nu >-1,
\end{eqnarray}
where 
\begin{equation}
\Phi(\xi, x, z)
=\prod_{a \in \xi \cap \{x\}^{\rm c}}
\frac{z-a}{x-a},
\label{eqn:entire}
\end{equation}
and ${\bf 1}(\omega)=1$ if $\omega$ is 
satisfied and ${\bf 1}(\omega)=0$ otherwise
(Proposition 2.1 in \cite{KT10} and 
Theorem 2.1 in \cite{KT11}).
The function $\Phi(\xi, x, z)$ is an entire function
of $z \in \C$ expressed by the Weierstrass canonical product
with genus 0, whose zeros are given by
$\supp \xi \cap \{x\}^{\rm c}$ \cite{Lev96,KT09}.

As direct consequences of Theorem \ref{Theorem:distributions},
we have the following equalities for
multitime correlation functions;
\begin{equation}
\bE_{N, \sigma^2}[
\rho^{\Xi}(t_1, \x^{(1)}_{N_1}; \dots;
t_M, \x^{(M)}_{N_M})]
=\rho^{N \delta_0}(t_1+\sigma^2, \x^{(1)}_{N_1};
\dots; t_M+\sigma^2, \x^{(M)}_{N_M})
\label{eqn:Th1b1}
\end{equation}
with $x^{(m)}_{N_m} \in \W_{N_m}^{\rm A}, 1 \leq m \leq M$,
and 
\begin{equation}
\bE^{(\nu)}_{N, \sigma^2}[
\rho_{\nu}^{\Xi}(t_1, \x^{(1)}_{N_1}; \dots;
t_M, \x^{(M)}_{N_M})]
=\rho_{\nu}^{N \delta_0}(t_1+\sigma^2, \x^{(1)}_{N_1};
\dots; t_M+\sigma^2, \x^{(M)}_{N_M})
\label{eqn:Th1b2}
\end{equation}
with $x^{(m)}_{N_m} \in \W_{N_m}^{+}, 1 \leq m \leq M$,
for any $M \in \N, 0 < t_1< \cdots < t_M < \infty$.
Therefore, infinite systems of equalities between
determinants of correlation kernels are obtained as
a corollary of Theorem \ref{Theorem:distributions}.

\begin{cor}
\label{Theorem:determinants}
For any $N, M \in \N, 0 < t_1 < \cdots < t_M < \infty,
N_m \in \{1,2, \dots, N\}, 1 \leq m \leq M$, $\sigma^2>0$,
\begin{eqnarray}
&& \bE_{N, \sigma^2} \left[
\det_{
\substack{1 \leq j \leq N_{m}, 1 \leq k \leq N_{n} \\ 1 \leq m, n \leq M}
}
[\mbK^{\Xi}(t_m, x_{j}^{(m)}; t_n, x_{k}^{(n)} )] \right]
\nonumber\\
&& \qquad 
=\det_{
\substack{1 \leq j \leq N_{m}, 1 \leq k \leq N_{n} \\ 1 \leq m, n \leq M}
}
[
\mbK^{N \delta_0}(t_m+\sigma^2, x_{j}^{(m)}; t_n+\sigma^2, x_{k}^{(n)} )],
\label{eqn:cor_eq1}
\end{eqnarray}
\begin{eqnarray}
&& \bE^{(\nu)}_{N, \sigma^2} \left[
\det_{
\substack{1 \leq j \leq N_{m}, 1 \leq k \leq N_{n} \\ 1 \leq m, n \leq M}
}
[\mbK^{\Xi}_{\nu}(t_m, x_{j}^{(m)}; t_n, x_{k}^{(n)} )] \right]
\nonumber\\
&& \qquad 
=\det_{
\substack{1 \leq j \leq N_{m}, 1 \leq k \leq N_{n} \\ 1 \leq m, n \leq M}
}
[
\mbK^{N \delta_0}_{\nu}
(t_m+\sigma^2, x_{j}^{(m)}; t_n+\sigma^2, x_{k}^{(n)} )].
\label{eqn:cor_eq2}
\end{eqnarray}
In particular, as $M=1, N_1=L \leq N$, 
\begin{eqnarray}
&& \bE_{N, \sigma^2} \left[
\det_{1 \leq j, k \leq L}
[\mbK^{\Xi}(t, x_{j}; t, x_{k})] \right]
=\det_{1 \leq j, k \leq L}
[\mbK^{N \delta_0}(t+\sigma^2, x_{j}; t+\sigma^2, x_{k})],
\label{eqn:cor_eq3}
\\
&& \bE^{(\nu)}_{N, \sigma^2} \left[
\det_{1 \leq j, k \leq L}
[\mbK^{\Xi}_{\nu}(t, x_{j}; t, x_{k})] \right]
=\det_{1 \leq j, k \leq L}
[\mbK_{\nu}^{N \delta_0}(t+\sigma^2, x_{j}; t+\sigma^2, x_{k})]
\label{eqn:cor_eq4}
\end{eqnarray}
hold for any $t >0, \sigma^2>0$.
\end{cor}
\vskip 0.5cm
\noindent
The proof of Theorem \ref{Theorem:distributions} is given in 
Sect. 3.1.

The main purpose of the present paper is to show that
the equalities in Corollary \ref{Theorem:determinants}
are nontrivial even in the special cases given by
(\ref{eqn:cor_eq3}) and (\ref{eqn:cor_eq4}), 
and from them the determinantal expressions
for the ensemble averages of $2n$-products of
characteristic polynomials of random matrices
are derived for any $n \in \N$.
We show two sets of determinantal expressions.
The first one is given by the following theorem.

\begin{thm}
\label{Theorem:Main1}
For any $N, n \in \N$, 
$\valpha=(\alpha_1, \alpha_2, \cdots, \alpha_{2n}) \in \C^{2n}$,
$\sigma^2>0$, 
\begin{eqnarray}
&&M_{\rm GUE}(2n, \valpha; N, \sigma^2)
=  
\frac{\gamma_{N,2n} \sigma^{n(2N+n)}}
{h_n(\alpha_1, \cdots, \alpha_n) h_n(\alpha_{n+1}, 
\cdots, \alpha_{2n})}
\nonumber\\
&& \qquad \times
\det_{1 \leq j, k \leq n}
\left[ \frac{1}{\alpha_j-\alpha_{n+k}}
\left| \begin{array}{cc}
H_{N+n} (\alpha_j/\sqrt{2 \sigma^2} ) 
& H_{N+n} ( \alpha_{n+k}/\sqrt{2 \sigma^2} ) \cr
H_{N+n-1} ( \alpha_j/\sqrt{2 \sigma^2} ) 
& H_{N+n-1}( \alpha_{n+k}/\sqrt{2 \sigma^2} )
\end{array} \right| \right] 
\label{eqn:chara1}
\end{eqnarray}
with
\begin{equation}
\gamma_{N,2n}=
2^{-n(2N+2n-1)/2} \prod_{\ell=2}^{n}
\frac{(N+n-\ell)!}{(N+n-1)!},
\label{eqn:gamma1}
\end{equation}
and for $\nu > -1$
\begin{eqnarray}
&& M^{(\nu)}(2n, \valpha; N, \sigma^2)
=\frac{\gamma_{N, 2n}^{(\nu)} (2 \sigma^2)^{n(2N+n)}}
{h_n(\alpha_1, \cdots, \alpha_n) h_n(\alpha_{n+1}, 
\cdots, \alpha_{2n})}
\nonumber\\
&& \quad \times
\det_{1 \leq j, k \leq n}
\left[ \frac{1}{\alpha_j-\alpha_{n+k}}
\left| \begin{array}{cc}
L^{\nu}_{N+n} \left(\alpha_j/2 \sigma^2 \right) 
& L^{\nu}_{N+n} \left( \alpha_{n+k}/2 \sigma^2 \right) \cr
L^{\nu}_{N+n-1} \left( \alpha_j/2 \sigma^2 \right) 
& L^{\nu}_{N+n-1}\left( \alpha_{n+k}/2 \sigma^2 \right)
\end{array} \right| \right],
\label{eqn:chara3}
\end{eqnarray}
with
\begin{equation}
\gamma_{N,2n}^{(\nu)}=(-1)^n 
\left( \frac{(N+n)!}{\Gamma(N+n+\nu)} \right)^{n-1}
\prod_{\ell=1}^{n-1} \Gamma(N+\nu+\ell)
\prod_{m=1}^{n+1} (N+m-1)!.
\label{eqn:gamma2}
\end{equation}
By setting $\nu \in \N_0, \nu=1/2$ and $\nu=-1/2$
in (\ref{eqn:chara3}),
the determinantal expressions are given
for $M^{(\nu)}_{\rm chGUE}, M_{\rm classC}$
and $M_{\rm classD}$ 
through (\ref{eqn:M2})-(\ref{eqn:M4}).
\end{thm}
\vskip 0.5cm
\noindent
Proof is given in Sect. 3.2.
The above expressions can be simplified by
using the following identity, which was 
given by Ishikawa et al.\cite{IOTZ06}.
For $n \geq 2, \x=(x_1, \dots, x_n), 
\y=(y_1, \dots, y_n),
\a=(a_1, \dots, a_n), \b=(b_1, \dots, b_n) \in \C^n$,
\begin{eqnarray}
&& \det_{1 \leq j, k \leq n}
\left[ \frac{1}{y_k-x_j} 
\left| \begin{array}{cc}
1 & a_j \cr 1 & b_k 
\end{array} \right| \right]
\nonumber\\
&=& \frac{(-1)^{n(n-1)/2}}
{\prod_{j=1}^{n} \prod_{k=1}^{n}(y_k-x_j)}
\left| \begin{array}{cccccccc}
1 & x_1 & \cdots & x_1^{n-1} & a_1 & a_1 x_1 & \cdots a_1 x_1^{n-1} \cr
1 & x_2 & \cdots & x_2^{n-1} & a_2 & a_2 x_2 & \cdots a_2 x_2^{n-1} \cr
  &     & \cdots &           &     &         & \cdots               \cr
1 & x_n & \cdots & x_n^{n-1} & a_n & a_n x_n & \cdots a_n x_n^{n-1} \cr
1 & y_1 & \cdots & y_1^{n-1} & b_1 & b_1 y_1 & \cdots b_1 y_1^{n-1} \cr
1 & y_2 & \cdots & y_2^{n-1} & b_2 & b_2 y_2 & \cdots b_2 y_2^{n-1} \cr
  &     & \cdots &           &     &         & \cdots               \cr
1 & y_n & \cdots & y_n^{n-1} & b_n & b_n y_n & \cdots b_n y_n^{n-1} \cr
\end{array} \right|.
\label{eqn:Ishikawa1}
\end{eqnarray}

\vskip 0.5cm
\noindent{\bf Remark 2.}
For $n \in \N, p,q \in \N_0$ satisfying $p+q=n$, and $\x, \a \in \C^n$,
denote by $V^{p,q}(\x;\a)$ the $n \times n$ matrix
with $j$-th row
$$
(1, x_j, \cdots, x_j^{p-1}, a_j, a_j x_j, \cdots,
a_j x_j^{q-1}).
$$
If $q=0$, then $p=n$ and 
$V^{n,0}(\x; \a)=V^{n,0}(\x)=(x_j^{k-1})_{1 \leq j, k \leq n}$
is the Vandermonde matrix and its determinant
$\det V^{n,0}(\x)$ is equal to the product
of differences of $n$ variables, $h_n(\x)$,
given by (\ref{eqn:Vandermonde}).
As a generalization of the Cauchy determinant
$$
\det_{1 \leq j, k \leq n}
\left( \frac{1}{x_j+y_k} \right)
= \frac{h_n(\x) h_n(\y)}
{\prod_{j=1}^{n} \prod_{k=1}^{n} (x_j+y_k)},
$$
$\x, \y \in \C^n$,
Ishikawa et al. \cite{IOTZ06} proved the
following equalities involving the generalized
Vandermonde determinants $V^{p,q}$.
Let $n \in \N$, $p,q \in \N_0$. 
For $\x=(x_1, \cdots, x_n), \y=(y_1, \cdots, y_n),
\a=(a_1, \cdots, a_n),\b=(b_1, \cdots, b_n) \in \C^n$
and
$\z=(z_1, \cdots, z_{p+q}), \c=(c_1, \cdots, c_{p+q}) \in \C^{p+q}$
\begin{eqnarray}
&& \det_{1 \leq j, k \leq n}
\left[ \frac{1}{y_k-x_j}
\det V^{p+1, q+1}(x_j, y_k, \z; a_j, b_k, \c) \right]
\nonumber\\
&=& \frac{(-1)^{n(n-1)/2}}
{\prod_{j=1}^{n} \prod_{k=1}^{n}(y_k-x_j)}
\det V^{p,q}(\z, \c)^{n-1}
\det V^{n+p, n+q}(\x, \y, \z; \a, \b, \c).
\label{eqn:Ishikawa2}
\end{eqnarray}
When $p=q=0$,
$\det V^{0,0}(\z,\c)=1$
and
$$
\det V^{1,1}(x_j, y_k; a_j, b_k)=
\left| \begin{array}{cc}
1 & a_j \cr 1 & b_k 
\end{array} \right|.
$$
Then as a special case of (\ref{eqn:Ishikawa2}),
(\ref{eqn:Ishikawa1}) is obtained.

For $\sigma^2>0$ define
\begin{eqnarray}
\label{eqn:hatH}
&& \widehat{H}_{\ell}(\alpha; \sigma^2)
=\left(\frac{\sigma^2}{2} \right)^{\ell/2}
H_{\ell} \left(\frac{\alpha}{\sqrt{2 \sigma^2}}\right),
\quad \alpha \in \R,
\\
\label{eqn:hatL}
&& \widehat{L}^{\nu}_{\ell}(\alpha; \sigma^2)
=(-2 \sigma^2)^{\ell} \ell !
L^{\nu}_{\ell} \left( \frac{\alpha}{2 \sigma^2} \right),
\quad \alpha \in \R_+,
\end{eqnarray}
$\ell \in \N_0$, which are both monic polynomials
of $\alpha$ with order $\ell$.
By using the identity (\ref{eqn:Ishikawa1})
and recurrence relations of Hermite polynomials
and Laguerre polynomials, we can prove
the following second set of determinantal expressions.

\begin{thm}
\label{Theorem:Main2}
For any $N, n \in \N, 
\valpha=(\alpha_1, \alpha_2, \cdots, \alpha_{2n}) \in \C^{2n}$,
$\sigma^2>0$, 
\begin{equation}
M_{\rm GUE}(2n, \valpha; N, \sigma^2)
= 
\frac{1}{h_{2n}(\valpha)}
\det_{1 \leq j, k \leq 2n}
\left[ \widehat{H}_{N+j-1}(\alpha_k; \sigma^2) \right],
\quad
\label{eqn:chara2}
\end{equation}
\begin{equation}
M_{\rm chGUE}^{(\nu)}(2n, \valpha; N, \sigma^2)
=
\frac{1}{h_{2n}(\valpha)}
\det_{1 \leq j, k \leq 2n}
\left[\widehat{L}^{\nu}_{N+j-1}(\alpha_k; \sigma^2) \right],
\quad \nu \in \N_0,
\label{eqn:chara4}
\end{equation}
\begin{eqnarray}
M_{\rm classC}(2n, \valpha; 2N, \sigma^2)
&=& \frac{1}{h_{2n}(\valpha^{\langle 2 \rangle})}
\det_{1 \leq j, k \leq 2n}
\left[ \widehat{L}^{1/2}_{N+j-1}
(\alpha_k^2; \sigma^2) \right]
\nonumber\\
\label{eqn:chara5}
&=& \frac{1}{h_{2n}(\valpha^{\langle 2 \rangle}) 
\prod_{j=1}^{2n} \alpha_j}
\det_{1 \leq j, k \leq 2n}
\left[ \widehat{H}_{2N+2j-1}
(\alpha_k; \sigma^2) \right],
\end{eqnarray}
and
\begin{eqnarray}
M_{\rm classD}(2n, \valpha; 2N, \sigma^2)
&=& \frac{1}{h_{2n}(\valpha^{\langle 2 \rangle})}
\det_{1 \leq j, k \leq 2n}
\left[ \widehat{L}^{-1/2}_{N+j-1}
(\alpha_k^2; \sigma^2) \right]
\nonumber\\
\label{eqn:chara6}
&=& \frac{1}{h_{2n}(\valpha^{\langle 2 \rangle})}
\det_{1 \leq j, k \leq 2n}
\left[ \widehat{H}_{2(N+j-1)}
(\alpha_k; \sigma^2) \right].
\end{eqnarray}
\end{thm}
\vskip 0.5cm

The determinantal expressions 
(\ref{eqn:chara2})-(\ref{eqn:chara6}) can be
obtained from the general formula given by
Br\'ezin and Hikami as Eq.(14) in \cite{BH00}.
(See also \cite{FW01,MN01,FS02,FS03} and
Sect.22.4 in \cite{Meh04}.)
Since our new expressions (\ref{eqn:chara1}) and (\ref{eqn:chara3})
are independently derived in the present paper,
if we combine the present result and that of
Br\'ezin and Hikami \cite{BH00},
the special case of identity (\ref{eqn:Ishikawa1})
of Ishikawa et al. \cite{IOTZ06} is concluded.

\vskip 0.5cm
\noindent{\bf Remark 3.}
In this paper, we derive the ensemble averages of
products of characteristic polynomials
of random matrices from the equalities 
(\ref{eqn:cor_eq3}) and (\ref{eqn:cor_eq4}).
As stated in Corollary \ref{Theorem:determinants},
these equalities are special cases
with $M=1$ of the systems of equalities
(\ref{eqn:cor_eq1}) and (\ref{eqn:cor_eq2}).
It will be an interesting problem to clarify
all the information involved in 
(\ref{eqn:cor_eq1}) and (\ref{eqn:cor_eq2})
(see \cite{Del10}).

\SSC{Proofs of Theorems}
\subsection{Proof of Theorem \ref{Theorem:distributions}}

By (\ref{eqn:Hermite1}) and the fact that $H_n(x/2)$ is 
a monic polynomial of order $n \in \N_0$ and thus
$$
\det_{1 \leq j, k \leq N}[H_{j-1}(x_k/2)]=h_N(\x),
$$
we obtain the equality
\begin{equation}
h^{(+)}_N(\sigma^2, \x; N \delta_0)
= \frac{\mu_{N, \sigma^2}(\xi)}{h_N(\xi)},
\quad \xi=\sum_{j=1}^{N} \delta_{x_j} \in \mM_{N,0},
\x \in \W_N^{\rm A}.
\label{eqn:h+}
\end{equation}
Combining this with (\ref{eqn:fhh1}) gives the equality
\begin{equation}
f(t_1, \x^{(1)}|\x) h^{(+)}_N(\sigma^2, \x; N \delta_0)
=h^{(+)}_N(t_1, \x^{(1)}; \xi) 
\mu_{N, \sigma^2}(\xi)
\label{eqn:fhhmu1}
\end{equation}
for any $t_1 >0, \xi=\sum_{j=1}^{N} \delta_{x_j} \in \mM_N,
\x^{(1)} \in \W_N^{\rm A}$.
Then (\ref{eqn:pxi}) with $\xi=N \delta_0$
gives
\begin{eqnarray}
&& p^{N \delta_0}(\sigma^2, \xi;
t_1+ \sigma^2, \xi^{(1)}; \dots; t_M+\sigma^2, \xi^{(M)})
\nonumber\\
&=& h^{(-)}_N(t_M+\sigma^2, \x^{(M)}; N \delta_0)
\prod_{m=1}^{M-1} f(t_{m+1}-t_{m}; \x^{(m+1)}|\x^{(m)})
f(t_1, \x^{(1)}|\x) h^{(+)}_N(\sigma^2, \x; N \delta_0)
\nonumber\\
&=& h^{(-)}_N(t_M, \x^{(M)}; \xi)
\prod_{m=1}^{M-1} f(t_{m+1}-t_{m}; \x^{(m+1)}|\x^{(m)})
h^{(+)}_N(t_1, \x^{(1)}; \xi)
\mu_{N, \sigma^2}(\xi)
\nonumber\\
&=& \mu_{N, \sigma^2}(\xi)
p^{\xi}(t_1, \xi^{(1)}; \dots; t_M, \xi^{(M)}),
\label{eqn:peq1}
\end{eqnarray}
where (\ref{eqn:h-}) was used.
Integrating the both sides with respect to
$\xi$ over $\mM_N$ according to (\ref{eqn:ENs}),
(\ref{eqn:Th1a1}) is obtained.
Similarly, the equality
\begin{equation}
f^{(\nu)}(t_1, \x^{(1)}|\x) h^{(\nu,+)}_N(\sigma^2, \x; N \delta_0)
=h^{(\nu,+)}_N(t_1, \x^{(1)}; \xi) 
\mu^{(\nu)}_{N, \sigma^2}(\xi)
\label{eqn:fhhmu2}
\end{equation}
is established for any $t_1 > 0,
\xi=\sum_{j=1}^{N} \delta_{x_j} \in \mM_N^+,
\xi^{(1)} \in \W_N^{+}$, and then we have
\begin{equation}
\mu^{(\nu)}_{N, \sigma^2}(\xi)
p_{\nu}^{\xi}(t_1, \xi^{(1)}; \dots; t_M, \xi^{(M)})
=p_{\nu}^{N \delta_0}(\sigma^2, \xi;
t_1+ \sigma^2, \xi^{(1)}; \dots; t_M+\sigma^2, \xi^{(M)}).
\label{eqn:peq2}
\end{equation}
Integrating the both sides with respect to
$\xi$ over $\mM_N^+$ according to (\ref{eqn:ENsn}),
(\ref{eqn:Th1a2}) is obtained. \qed

\subsection{Proof of Theorem \ref{Theorem:Main1}}

First we derive the expression (\ref{eqn:chara1})
from (\ref{eqn:cor_eq3}).
We observe that the LHS of (\ref{eqn:cor_eq3}) is
written as follows.

\begin{lem}
\label{Theorem:Lemma3_1}
For any $\sigma^2>0, t >0, 1 \leq L \leq N,
\x_{L}=(x_1, \dots, x_L) \in \W_{L}^{\rm A}$,
\begin{eqnarray}
&& \bE_{N, \sigma^2}
\left[ \det_{1 \leq j, k \leq L}
[\mbK^{\Xi}(t,x_j; t, x_k)] \right]
\nonumber\\
&=& \int_{\R^{L}} d \u \int_{\R^{L}} d \w \,
e^{-|\u|^2/2 \sigma^2}
\prod_{j=1}^{L} p(t, x_j|u_j)
p(-t, \im w_j|x_j)
\frac{
\sigma^{-L(2N-L)}}
{(2 \pi)^{L/2} \prod_{\ell=1}^{L}(N-\ell)!}
\nonumber\\
&& \qquad \times
h_{L}(\u) h_{L}(\im \w)
M_{\rm GUE}(2L, (u_1, \dots, u_{L}, \im w_1, \dots, \im w_{L});
N-L, \sigma^2).
\label{eqn:eqLem3_1}
\end{eqnarray}
\end{lem}
\noindent{\it Proof} \quad
Assume that $\xi=\sum_{j=1}^N \delta_{a_j}$ with
$\a=(a_1, \dots, a_N) \in \W_N^{\rm A}$.
Then
\begin{eqnarray}
&& \rho(t, \x_L)= \det_{1 \leq j, k \leq L} 
[\mbK^{\xi}(t, x_j; t, x_k)] 
\nonumber\\
&& 
=\det_{1 \leq j, k \leq L} \left[
\sum_{q=1}^{N} p(t,x_j|a_q) \int_{\R} dw \,
p(-t, \im w|x_k)
\prod_{1 \leq \ell \leq N, \ell \not=q}
\frac{\im w-a_{\ell}}{a_q-a_{\ell}} \right]
\nonumber\\
&& 
= L! \sum_{1 \leq q_1 < \cdots < q_N \leq N}
\int_{\R^L} d \w \,
\prod_{j=1}^{L} \left\{
p(-t, \im w_j|x_j)
\prod_{1 \leq \ell_j \leq N, \ell_j \not= q_j}
\frac{\im w_j-a_{\ell_j}}{a_{q_j}-a_{\ell_j}} \right\}
\det_{1 \leq j, k \leq L} 
[ p(t, x_j|a_{q_k})],
\nonumber
\end{eqnarray}
where multilinearlity and antisymmetry property of determinant
have been used.
Let $\I_N=\{1,2,\dots,N\}, \I_L=\{1,2, \dots, L\}$.
For a given ordered set of indices
$\q_L=\{q_1, q_2, \dots, q_L\},
1 \leq q_1 < \cdots < q_L \leq N$,
we see
\begin{eqnarray}
&& h_N(\a)^2 
\prod_{j=1}^L \prod_{1 \leq \ell_j \leq N, 
\ell_j \not= q_j}
\frac{\im w_j-a_{\ell_j}}{a_{q_j}-a_{\ell_j}}
\nonumber\\
&& \qquad
= (-1)^{L(L-1)/2} h_{N-L}((a_j)_{j \in \I_N \setminus \q_L})
\nonumber\\
&& \qquad \quad
\times \prod_{k \in \I_L} \prod_{j \in \I_N \setminus \q_L}
(a_{q_k}-a_j)(\im w_k-a_j)
\times \prod_{k \in \I_L} \prod_{j \in \I_L, j \not=k}
(\im w_k-a_{p_j}).
\nonumber
\end{eqnarray}
We also note that
$$
\frac{\sigma^{-N^2}}{N! C_N}
=\frac{\sigma^{-(N-L)^2}}{(N-L)! C_{N-L}}
\times \frac{\sigma^{-L(2N-N')}}
{(2\pi)^{L/2} \prod_{n=0}^{L-1} (N-n)!}.
$$
Then the LHS of (\ref{eqn:eqLem3_1}) becomes
\begin{eqnarray}
&& 
\int_{\R^N} d \a \,
\frac{\sigma^{-N^2}}{N! C_N} e^{-|\a|^2/2\sigma}
h_N(\a)^2 \rho^{\xi}(t, \x_L)
\nonumber\\
&=& \frac{(-1)^{L(L-1)/2} \sigma^{-L(2N-L)}}
{(2\pi)^{L/2} \prod_{n=0}^{L-1} (N-n)!}
L! \sum_{1 \leq q_1 < \cdots < q_L \leq N}
\int_{\R^L} d \w
\prod_{j=1}^{L} p(-t, \im w_j|x_j)
\prod_{k=1}^{L} \int_{\R} da_{q_k} \,
e^{-a_{q_k}^2/2 \sigma^2}
\nonumber\\
&\times& \prod_{\ell \in \I_L} \prod_{m \in \I_L, m \not=\ell}
(\im w_{\ell}-a_{q_m})
\det_{1 \leq j, k \leq L} [p(t, x_j|a_{q_k})]
\bE_{N-L, \sigma^2} \left[
\prod_{k \in \I_L} \prod_{j \in \I_L \setminus \q_L}
(a_{q_k}-a_j)(\im w_k-a_j) \right]
\nonumber\\
&=& 
\frac{(-1)^{L(L-1)/2} \sigma^{-L(2N-L)}}
{(2\pi)^{L/2} \prod_{n=1}^L (N-n)!}
\int_{\R^L} d \w \,
\prod_{j=1}^{L} p(-t, \im w_j|x_j)
\int_{\R^L} d \v \, e^{-|\v|^2/2 \sigma^2}
\prod_{k=1}^{L} \prod_{1 \leq j \leq L, j \not= k}
(\im w_k-v_j)
\nonumber\\
&\times&
\det_{1 \leq j, k \leq L}[p(t, x_j|v_k)]
M_{\rm GUE}(2L, (v_1, \dots, v_L, \im w_1, \dots, \im w_L);
N-L, \sigma^2),
\nonumber
\end{eqnarray}
where we have replaced the integral variables 
$(a_{q_1}, \dots, a_{q_L})$ by $(v_1, \dots, v_L) \equiv \v$.
By definition
$$
\det_{1 \leq j, k \leq L}
[p(t, x_j|v_k)]
=\sum_{\tau \in \S_L} {\rm sgn}(\tau)
\prod_{j=1}^L p(t, x_j|v_{\tau(j)}),
$$
where $\S_L$ denotes the collection of all permutation
of $(1,2, \dots, L)$.
For each $\tau \in \S_L$, set
$v_{\tau(j)}=u_j, 1 \leq j \leq L$, that is,
$v_j=u_{\tau^{-1}(j)}, 1 \leq j \leq L$.
Then the above equals
\begin{eqnarray}
&& \frac{(-1)^{L(L-1)/2} \sigma^{-L(2N-L)}}
{(2\pi)^{L/2} \prod_{n=1}^L (N-n)!}
\int_{\R^L} d \u 
\int_{\R^L} d \w \, 
e^{-|\u|^2/2 \sigma^2}
\prod_{j=1}^{L} p(t, x_j|u_j) p(-t, \im w_j|x_j)
\nonumber\\
&& \quad \times
\sum_{\tau \in \S_L} {\rm sgn}(\tau)
\prod_{k=1}^{L} \prod_{1 \leq j \leq L, j \not=k}
(\im w_k-u_{\tau^{-1}(j)})
\nonumber\\
&& \qquad \quad \times
M_{\rm GUE}(2L, (u_{\tau^{-1}(1)}, \dots,
u_{\tau^{-1}(L)}, \im w_1, \dots, \im w_L);
N-L, \sigma^2).
\nonumber
\end{eqnarray}
By definition (\ref{eqn:PGUE1}),
$M_{\rm GUE}(m, \valpha; N, \sigma^2)$ is symmetric in $\valpha$,
and thus
\begin{eqnarray}
&&
M_{\rm GUE}(2L, (u_{\tau^{-1}(1)}, \dots,
u_{\tau^{-1}(L)}, \im w_1, \dots, \im w_L);
N-L, \sigma^2)
\nonumber\\
&& \quad =
M_{\rm GUE}(2L, (u_1, \dots, u_L, \im w_1, \dots, \im w_L);
N-L, \sigma^2).
\nonumber
\end{eqnarray}
We can confirm that for any $n \geq 2, \x, \y \in \R^n$
$$
\sum_{\kappa \in \S_n} {\rm sgn}(\kappa)
\prod_{k=1}^n \prod_{1 \leq j \leq n, j \not=k}
(x_k-y_{\kappa(j)})
=(-1)^{[n/2]} h_n(\x) h_n(\y).
$$
We can see that 
$[L/2]+L(L-1)/2$ is even for any $L \in \N$.
Then (\ref{eqn:eqLem3_1}) is obtained. \qed
\vskip 0.5cm

Then we consider the RHS of (\ref{eqn:cor_eq3}).
By the first expression of (\ref{eqn:Kxi}) and
the fact (\ref{eqn:Hermite1}), we have
the Hermite kernel \cite{Meh04,For10}
\begin{equation}
\mbK^{N \delta_0}(t+\sigma^2,x; t+\sigma^2,y)
=\frac{e^{-x^2/2(t+\sigma^2)}}
{\sqrt{2\pi (t+\sigma^2)}}
\sum_{n=0}^{N-1}
\frac{1}{2^n n!} H_n \left( \frac{x}{\sqrt{2(t+\sigma^2)}} \right)
H_n \left( \frac{y}{\sqrt{2(t+\sigma^2)}} \right),
\label{eqn:KN0a}
\end{equation}
$t>0, \sigma^2>0, x,y \in \R$.
By the extended version of Chapman-Kolmogorov equation
(\ref{eqn:CK1}), the following integral formulas are derived.

\begin{lem}
\label{Theorem:Lemma3_2}
For $t >0, \sigma^2>0, n \in \N_0, x, y \in \R$
\begin{eqnarray}
&& \int_{\R} du \, p(t, x|u) 
H_n \left( \frac{u}{\sqrt{2 \sigma^2}} \right)
e^{-u^2/2\sigma^2}
=\left(\frac{\sigma^2}{t+\sigma^2}\right)^{(n+1)/2}
H_n \left(\frac{x}{\sqrt{2(t+\sigma^2)}}\right)
e^{-x^2/2(t+\sigma^2)},
\nonumber\\
\label{eqn:IntA1}
\\
\label{eqn:IntA2}
&& \int_{\R} dv \,
H_n \left( \frac{\im v}{\sqrt{2 \sigma^2}} \right)
p(-t, \im v|x)
=\left(\frac{t+\sigma^2}{\sigma^2} \right)^{n/2}
H_n \left( \frac{x}{\sqrt{2(t+\sigma^2)}}\right),
\\
&& \mbK^{N \delta_0}(t+\sigma^2,x; t+\sigma^2,y)
= \frac{1}{2^N (N-1)! \sqrt{\pi}}
\int_{\R} du \int_{\R} dv \,
e^{-u^2/2\sigma^2} p(t,x|u) p(-t,\im v|y)
\nonumber\\
\label{eqn:IntA3}
&& \qquad \qquad \qquad \qquad \qquad \qquad \qquad
\times \frac{1}{u-\im v}
\left| \begin{array}{ll}
H_N(u/\sqrt{2\sigma^2}) & H_N(\im v/\sqrt{2\sigma^2}) \cr
H_{N-1}(u/\sqrt{2\sigma^2}) & H_{N-1}(\im v/\sqrt{2\sigma^2})
\end{array} \right|.
\end{eqnarray}
\end{lem}
\noindent{\it Proof} \quad
Since the equality
$$
\int_{\R} du \, p(t, x|u) p(\sigma^2, u|s)
=p(t+\sigma^2, x|s)
$$
holds for $t>0, \sigma^2>0, x \in \R, s \in \C$,
(\ref{eqn:phi1a}) gives
$$
\int_{R}du \, p(t,x|u) \phi_n^{(+)}(\sigma^2,u;\xi)
=\phi_n^{(+)}(t+\sigma^2, x; \xi),
$$
$n \in \N_0, t >0, \sigma^2>0, x \in \R$.
If we set $\xi=N \delta_0$ and use the fact (\ref{eqn:Hermite1}),
we obtain (\ref{eqn:IntA1}).
Similarly, by the equality
\begin{eqnarray}
&& \int_{\R} dv \,
p(-\sigma^2, \im s| \im v) p(-t, \im v|x)
=\int_{\R} dv \, 
p(\sigma^2, s|v) p(t,v|-\im x)
\nonumber\\
&& \qquad = p(t+\sigma^2, s|-\im x)
=p(-(t+\sigma^2), \im s|x),
\nonumber
\end{eqnarray}
$t>0, \sigma^2>0, x \in \R$, and by the definition
(\ref{eqn:phi1b}), we see
$$
\int_{\R}dv \, \phi_n^{(-)}(\sigma^2, \im v; \xi)
p(-t, \im v|x)
= \phi_n^{(-)}(t+\sigma^2,x; \xi).
$$
If we set $\xi=N \delta_0$ and use the fact (\ref{eqn:Hermite1}),
we obtain (\ref{eqn:IntA2}).
Inserting (\ref{eqn:IntA1}) and (\ref{eqn:IntA2}) into
the RHS of (\ref{eqn:KN0a}), we have
\begin{eqnarray}
\mbK^{N \delta_0}(t+\sigma^2, x; t+\sigma^2, t)
&=&\frac{1}{\sqrt{2 \pi \sigma^2}}
\int_{\R} du \int_{\R} dv \,
e^{-u^2/2\sigma^2} p(t, x|u) p(-t, \im v|y)
\nonumber\\
&& \quad
\times \sum_{n=0}^{N-1} \frac{1}{2^n n!} 
H_n \left( \frac{u}{\sqrt{2 \sigma^2}}\right)
H_{n} \left( \frac{\im v}{\sqrt{2 \sigma^2}} \right).
\nonumber
\end{eqnarray}
Here we use the Christoffel-Darboux formula 
(see, for example, \cite{AAR99}), 
$$
\sum_{n=0}^{N-1} \frac{1}{2^n n!} H_n(x) H_n(y)
=\frac{1}{2^N (N-1)! (x-y)} 
\left| \begin{array}{ll}
H_{N}(x) & H_N(y) \cr
H_{N-1}(x) & H_{N-1}(y) 
\end{array} \right|,
$$
and then (\ref{eqn:IntA3}) is obtained. \qed
\vskip 0.5cm

The RHS of (\ref{eqn:cor_eq3}) is thus written as follows.

\begin{lem}
\label{Theorem:Lemma3_3}
For any $\sigma^2>0, t >0, 1 \leq L \leq N,
\x_{L}=(x_1, \dots, x_L) \in \W_{L}^{\rm A}$,
\begin{eqnarray}
&& \det_{1 \leq j, k \leq L}
[\mbK^{N \delta_0}(t+\sigma^2, x_j; t+\sigma^2, x_k)]
\nonumber\\
&& \quad = \int_{\R^L} d \u \int_{\R^L} d \w \,
e^{-|\u|^2/2\sigma^2}
\prod_{j=1}^{L} p(t, x_j|u_j) p(-t, \im w_j|x_j)
\nonumber\\
&& \qquad 
\times 
\frac{1}{(2^N (N-1)! \sqrt{\pi})^L}
\det_{1 \leq j, k \leq L} \left[
\frac{1}{u_j-\im w_k} 
\left| \begin{array}{ll}
H_N(u_j/\sqrt{2 \sigma^2}) & H_{N}(\im w_k/\sqrt{2 \sigma^2}) \cr
H_{N-1}(u_j/\sqrt{2 \sigma^2}) & 
H_{N-1}(\im w_k/\sqrt{2 \sigma^2}) 
\end{array} \right| \right].
\nonumber\\
\label{eqn:eqLem3_3}
\end{eqnarray}
\end{lem}
\noindent{\it Proof} \quad
By definition of determinant, the LHS of (\ref{eqn:eqLem3_3})
is given by
\begin{eqnarray}
&& \frac{1}{(2^N (N-1)! \sqrt{\pi})^L}
\sum_{\tau \in \S_L} {\rm sgn}(\tau)
\int_{\R^L} d \u \int_{\R^L} d \v \,
e^{-|\u|^2/2 \sigma^2} 
\nonumber\\
&& \quad \times
\prod_{j=1}^{L} \left\{
p(t,x_j|u_j) \frac{1}{u_j-\im v_j}
\left| \begin{array}{ll}
H_N(u_j/\sqrt{2 \sigma^2}) & H_N(\im v_j/\sqrt{2 \sigma^2}) \cr
H_{N-1}(u_j/\sqrt{2 \sigma^2}) & 
H_{N-1}(\im v_j/\sqrt{2 \sigma^2})
\end{array} \right|
p(-t, \im v_j|x_{\tau(j)}) \right\}.
\nonumber
\end{eqnarray}
Note that for each permutation $\tau \in \S_L$
$$
\prod_{j=1}^L p(-t, \im v_j|x_{\tau(j)})
=\prod_{j=1}^{L} p(-t, \im v_{\tau^{-1}(j)}|x_j).
$$
We set $v_{\tau^{-1}(j)}=w_j, 1 \leq j \leq N$,
that is, $v_j=w_{\tau(j)}, 1 \leq j \leq N$.
Then the above equals
\begin{eqnarray}
&& \frac{1}{2^{NL}((N-1)!)^L \pi^{L/2}}
\int_{\R^L} d \u \int_{\R^L} d \w \,
e^{-|\u|^2/2 \sigma^2} 
\prod_{j=1}^{L} \{
p(t,x_j|u_j) p(-t, \im w_j|x_j) \}
\nonumber\\
&& \quad \times
\sum_{\tau \in \S_L} {\rm sgn}(\tau)
\prod_{j=1}^{L} \left\{
\frac{1}{u_j-\im w_{\tau(j)}}
\left| \begin{array}{ll}
H_N(u_j/\sqrt{2 \sigma^2}) & H_N(\im w_{\tau(j)}/\sqrt{2 \sigma^2}) \cr
H_{N-1}(u_j/\sqrt{2 \sigma^2}) & 
H_{N-1}(\im w_{\tau(j)}/\sqrt{2 \sigma^2})
\end{array} \right| \right\},
\nonumber
\end{eqnarray}
which is the RHS of (\ref{eqn:eqLem3_3}). \qed
\vskip 0.5cm

\noindent{\it 
Proof of Eq.(\ref{eqn:chara1}) in Theorem \ref{Theorem:Main1}
} \quad
Since the equality (\ref{eqn:cor_eq3}) 
in Corollary \ref{Theorem:determinants},
holds for any $t >0, \sigma^2 >0$, 
Lemmas \ref{Theorem:Lemma3_1} and \ref{Theorem:Lemma3_3},
imply that 
the integrand of multiple Gaussian integral 
in the RHS of (\ref{eqn:eqLem3_1})
is equal to that in the RHS of (\ref{eqn:eqLem3_3}).
By replacing $L$ and $N-L$ by
$n$ and $N$, respectively, and
$(u_1, \dots, u_L, \im w_1, \dots, \im w_L)$
by $(\alpha_1, \dots, \alpha_n, \alpha_{n+1},
\dots, \alpha_{2n})=\alpha \in \C^{2n}$,
we obtain (\ref{eqn:chara1}). \qed
\vskip 0.5cm

Next we prove (\ref{eqn:chara3}). 
The LHS of (\ref{eqn:cor_eq4}) is
written as follows.

\begin{lem}
\label{Theorem:Lemma3_4}
For any $\sigma^2>0, t >0, 1 \leq L \leq N,
\x_{L}=(x_1, \dots, x_L) \in \W_{L}^{+}$,
\begin{eqnarray}
&& \bE^{(\nu)}_{N, \sigma^2}
\left[ \det_{1 \leq j, k \leq L}
[\mbK_{\nu}^{\Xi}(t,x_j; t, x_k)] \right]
\nonumber\\
&=& \int_{\R_+^{L}} d \u \int_{\R_-^{L}} d \w \,
\prod_{j=1}^{L} 
u_{j}^{\nu} e^{-u_{j}/2 \sigma^2}
p^{(\nu)}(t, x_j|u_j)
p^{(\nu)}(-t, w_j|x_j)
\frac{(2 \sigma^2)^{-L(2N-L+\nu)}}
{\prod_{\ell=1}^L \Gamma(N+\nu+1-\ell) (N-\ell)!}
\nonumber\\
&& \qquad \times h_{L}(\u) h_{L}(\w)
M^{(\nu)}(2L, (u_1, \dots, u_{L}, w_1, \dots, w_{L});
N-L, \sigma^2).
\label{eqn:eqLem3_4}
\end{eqnarray}
\end{lem}
\vskip 0.5cm
\noindent
Since we can prove this lemma in the similar way to
Lemma \ref{Theorem:Lemma3_1}, we omit the proof.

Then we consider the RHS of (\ref{eqn:cor_eq4}).
By the first expression of (\ref{eqn:Kxinu}) and
the fact (\ref{eqn:Laguerre1}), we have
the Laguerre kernel \cite{For10}
\begin{eqnarray}
&& \mbK_{\nu}^{N \delta_0}(t+\sigma^2,x; t+\sigma^2,y)
\nonumber\\
&& \qquad 
=\frac{x^{\nu} e^{-x/2(t+\sigma^2)}}
{\{2(t+\sigma^2)\}^{\nu+1}}
\sum_{n=0}^{N-1}
\frac{n!}{\Gamma(n+\nu+1)} 
L^{\nu}_n \left( \frac{x}{2(t+\sigma^2)} \right)
L^{\nu}_n \left( \frac{y}{2(t+\sigma^2)} \right),
\label{eqn:KN0nua}
\end{eqnarray}
$t>0, \sigma^2>0, x,y \in \R_+$.
By the extended version of Chapman-Kolmogorov equation
(\ref{eqn:CK2}), the following integral formulas are derived.

\begin{lem}
\label{Theorem:Lemma3_5}
For $t >0, \sigma^2>0, n \in \N_0, x, y \in \R_+$
\begin{eqnarray}
&& \int_{\R_+} du \, 
p^{(\nu)}(t, x|u) 
L^{\nu}_n \left( \frac{u}{2 \sigma^2} \right)
u^{\nu} e^{-u/2 \sigma^2}
= \left(\frac{\sigma^2}{t+\sigma^2}\right)^{n+\nu+1}
L^{\nu}_n \left(\frac{x}{2(t+\sigma^2)}\right)
x^{\nu}e^{-x/2(t+\sigma^2)},
\nonumber\\
\label{eqn:IntB1}
\\
\label{eqn:IntB2}
&& \int_{\R_-} dv \,
L^{\nu}_n \left( \frac{v}{2 \sigma^2} \right)
p^{(\nu)}(-t, v|x)
= \left(\frac{t+\sigma^2}{\sigma^2} \right)^{n}
L^{\nu}_n \left( \frac{x}{2(t+\sigma^2)}\right),
\\
&& \mbK_{\nu}^{N \delta_0}(t+\sigma^2,x; t+\sigma^2,y)
\nonumber\\
&& \qquad 
= -\frac{N!}{(2 \sigma^2)^{\nu} \Gamma(N+\nu)}
\int_{\R_+} du \int_{\R_-} dv \,
u^{\nu} e^{-u/2\sigma^2} p^{(\nu)}(t, x|u) p^{(\nu)}(-t, v|y)
\nonumber\\
\label{eqn:IntB3}
&& \qquad \qquad \qquad \qquad \qquad \qquad
\times \frac{1}{u-v}
\left| \begin{array}{ll}
L^{\nu}_N(u/2\sigma^2) & L^{\nu}_N(v/2\sigma^2) \cr
L^{\nu}_{N-1}(u/2\sigma^2) & L^{\nu}_{N-1}(v/2\sigma^2)
\end{array} \right|.
\end{eqnarray}
\end{lem}
\noindent{\it Proof} \quad
Since the equality
$$
\int_{\R_+} du \, p^{(\nu)}(t, x|u) p^{(\nu)}(\sigma^2, u|s)
=p^{(\nu)}(t+\sigma^2, x|s)
$$
holds for $t>0, \sigma^2>0, x \in \R_+, s \in \C$,
(\ref{eqn:phinu1a}) gives
$$
\int_{R_+}du \, p^{(\nu)}(t,x|u) \phi_n^{(\nu,+)}(\sigma^2,u;\xi)
=\phi_n^{(\nu, +)}(t+\sigma^2, x; \xi),
$$
$n \in \N_0, t >0, \sigma^2>0, x \in \R_+$.
If we set $\xi=N \delta_0$ and use (\ref{eqn:Laguerre1}),
we obtain (\ref{eqn:IntB1}).
Similarly, by the equality
\begin{eqnarray}
&& \int_{\R_-} dv \,
p^{(\nu)}(-\sigma^2, -w|v) p^{(\nu)}(-t, v|x)
= \int_{\R_+} du \, p^{(\nu)}(\sigma^2, w|u)
p^{(\nu)}(t, u|-x)
\nonumber\\
&& \qquad
=p^{(\nu)}(t+\sigma^2, w|-x)
=p^{(\nu)}(-(t+\sigma^2), -w|x),
\nonumber
\end{eqnarray}
$t>0, \sigma^2>0, x, w \in \R_+$, and by the definition
(\ref{eqn:phinu1b}), we see
$$
\int_{\R_-}dv \, \phi_n^{(\nu,-)}(\sigma^2, v; \xi)
p^{(\nu)}(-t, v|x)
= \phi_n^{(\nu,-)}(t+\sigma^2,x; \xi).
$$
If we set $\xi=N \delta_0$ and use (\ref{eqn:Laguerre1}),
we obtain (\ref{eqn:IntB2}).
Inserting (\ref{eqn:IntB1}) and (\ref{eqn:IntB2}) into
the RHS of (\ref{eqn:KN0nua}), we have
\begin{eqnarray}
\mbK_{\nu}^{N \delta_0}(t+\sigma^2, x; t+\sigma^2, t)
&=&\frac{1}{(2 \sigma^2)^{\nu+1}}
\int_{\R_+} du \int_{\R_-} dv \,
u^{\nu} e^{-u/2\sigma^2} p^{(\nu)}(t, x|u) p^{(\nu)}(-t, v|y)
\nonumber\\
&& \quad
\times \sum_{n=0}^{N-1} \frac{n!}{\Gamma(n+\nu+1)} 
L^{\nu}_n \left( \frac{u}{2 \sigma^2}\right)
L^{\nu}_{n} \left( \frac{v}{2 \sigma^2} \right).
\nonumber
\end{eqnarray}
Here we use the Christoffel-Darboux formula \cite{AAR99}
$$
\sum_{n=0}^{N-1} \frac{n!}{\Gamma(n+\nu+1)} L^{\nu}_n(x) L^{\nu}_n(y)
=-\frac{N!}{\Gamma(N+\nu) (x-y)} 
\left| \begin{array}{ll}
L^{\nu}_{N}(x) & L^{\nu}_N(y) \cr
L^{\nu}_{N-1}(x) & L^{\nu}_{N-1}(y) 
\end{array} \right|.
$$
Then (\ref{eqn:IntB3}) is obtained. \qed
\vskip 0.5cm

The RHS of (\ref{eqn:cor_eq4}) is thus written as follows.

\begin{lem}
\label{Theorem:Lemma3_6}
For any $\sigma^2>0, t >0, 1 \leq L \leq N,
\x_{L}=(x_1, \dots, x_L) \in \W_{L}^{+}$,
\begin{eqnarray}
&& \det_{1 \leq j, k \leq L}
[\mbK_{\nu}^{N \delta_0}(t+\sigma^2, x_j; t+\sigma^2, x_k)]
\nonumber\\
&& \quad = \int_{\R_+^L} d \u \int_{\R_-^L} d \w \,
\prod_{j=1}^{L} u_j^{\nu} e^{-u_j/2 \sigma^2} p^{(\nu)}(t, x_j|u_j)
p^{(\nu)}(-t, w_j|x_j)
\nonumber\\
&& \qquad 
\times \left(
-\frac{N!}{(2 \sigma^2)^{\nu} \Gamma(N+\nu)} \right)^L
\det_{1 \leq j, k \leq L} \left[
\frac{1}{u_j-w_k} 
\left| \begin{array}{ll}
L^{\nu}_N(u_j/2 \sigma^2) & L^{\nu}_{N}(w_k/2 \sigma^2) \cr
L^{\nu}_{N-1}(u_j/2 \sigma^2) & 
L^{\nu}_{N-1}(w_k/2 \sigma^2) 
\end{array} \right| \right].
\nonumber\\
\label{eqn:eqLem3_6}
\end{eqnarray}
\end{lem}
\vskip 0.5cm

\noindent{\it 
Proof of Eq.(\ref{eqn:chara3}) in Theorem \ref{Theorem:Main1}
} \quad
Since the equality (\ref{eqn:cor_eq4}) in Corollary 
\ref{Theorem:determinants},
holds for any $t >0, \sigma^2 >0$, 
Lemmas \ref{Theorem:Lemma3_4} and \ref{Theorem:Lemma3_6},
imply that
the integrand of multiple Gaussian integral 
in the RHS of (\ref{eqn:eqLem3_4})
is equal to that in the RHS of (\ref{eqn:eqLem3_6}).
By replacing $L$ and $N-L$ by
$n$ and $N$, respectively, and
$(u_1, \dots, u_L, w_1, \dots, w_L)$
by $(\alpha_1, \dots, \alpha_n, \alpha_{n+1},
\dots, \alpha_{2n})=\alpha \in \C^{2n}$,
we obtain (\ref{eqn:chara3}). \qed

\subsection{Proof of Theorem \ref{Theorem:Main2}}

First we prove that (\ref{eqn:chara2}) is equal to 
(\ref{eqn:chara1}).
Since
\begin{eqnarray}
&& \frac{1}{\alpha_j-\alpha_{n+k}}
\left| \begin{array}{ll}
H_{N+n}(\alpha_j/\sqrt{2 \sigma^2}) &
H_{N+n}(\alpha_{n+k}/\sqrt{2 \sigma^2}) \cr
H_{N+n-1}(\alpha_j/\sqrt{2 \sigma^2}) &
H_{N+n-1}(\alpha_{n+k}/\sqrt{2 \sigma^2}) 
\end{array} \right|
\nonumber\\
&& \quad =
H_{N+n-1} \left( \frac{\alpha_j}{\sqrt{2 \sigma^2}} \right)
H_{N+n-1} \left( \frac{\alpha_{n+k}}{\sqrt{2 \sigma^2}}\right)
\nonumber\\
&& \qquad \quad \times
\frac{1}{\alpha_{n+k}-\alpha_j}
\left| \begin{array}{ll}
1 & H_{N+n}(\alpha_j/\sqrt{2 \sigma^2})/H_{N+n-1}(\alpha_j/\sqrt{2 \sigma^2})
\cr
1 & H_{N+n}(\alpha_{n+k}/\sqrt{2 \sigma^2})
/H_{N+n-1}(\alpha_{n+k}/\sqrt{2 \sigma^2})
\end{array} \right|,
\nonumber
\end{eqnarray}
we can apply the identity (\ref{eqn:Ishikawa1})
and obtain the following,
\begin{eqnarray}
&& \det_{1 \leq j, k \leq n}
\left[ \frac{1}{\alpha_j-\alpha_{n+k}}
\left| \begin{array}{ll}
H_{N+n}(\alpha_j/\sqrt{2 \sigma^2}) &
H_{N+n}(\alpha_{n+k}/\sqrt{2 \sigma^2}) \cr
H_{N+n-1}(\alpha_j/\sqrt{2 \sigma^2}) &
H_{N+n-1}(\alpha_{n+k}/\sqrt{2 \sigma^2}) 
\end{array} \right| \right]
\nonumber\\
&=& \prod_{j=1}^n \left\{
H_{N+n-1} \left( \frac{\alpha_j}{\sqrt{2 \sigma^2}} \right)
H_{N+n-1} \left( \frac{\alpha_{n+k}}{\sqrt{2 \sigma^2}}\right) \right\}
\frac{(-1)^{n(n-1)/2}}
{\prod_{j=1}^{n} \prod_{k=1}^{n}(\alpha_{n+k}-\alpha_{n})}
\nonumber\\
&& \times
V^{n,n} \left(
(\alpha_1, \dots, \alpha_n),
(\alpha_{n+1}, \dots, \alpha_{2n});
\left( \frac{H_{N+n}(\frac{\alpha_j}{\sqrt{2 \sigma^2}})}
{H_{N+n-1}(\frac{\alpha_j}{\sqrt{2 \sigma^2}})} \right)_{j=1}^n,
\left( \frac{H_{N+n}(\frac{\alpha_{n+j}}{\sqrt{2 \sigma^2}})}
{H_{N+n-1}(\frac{\alpha_{n+j}}{\sqrt{2 \sigma^2}})} 
\right)_{j=1}^n \right)
\nonumber\\
&=& 
\frac{(-1)^{n(n-1)/2}}
{\prod_{j=1}^{n} \prod_{k=1}^{n}(\alpha_{n+k}-\alpha_{n})}
\det [A_{N+n-1} \quad A_{N+n} ],
\nonumber
\end{eqnarray}
where $A_{\ell}$ is the $2n \times n$ rectangular
matrix whose $(j,k)$-element is
given by $\alpha_j^{k-1} H_{\ell}
(\alpha_j/\sqrt{2 \sigma^2})$,
$1 \leq j \leq 2n, 1 \leq k \leq n$.
Here we use the recurrence relation of Hermite
polynomials
\begin{equation}
H_{n+1}(x)=2x H_n(x)-2n H_{n-1}(x).
\label{eqn:rec1}
\end{equation}
Then we find that
\begin{eqnarray}
&& \alpha^{k-1} H_{N+n-1}\left( \frac{\alpha}{\sqrt{2 \sigma^2}} \right)
=
(2 \sigma^2)^{(k-1)/2} \prod_{\ell=1}^{k-1} (N+n-\ell) H_{N+n-k}
\left( \frac{\alpha}{\sqrt{2 \sigma^2}} \right)
\nonumber\\
&& \quad + \mbox{linear combination of
$\left\{H_{\ell}
(\alpha/\sqrt{2 \sigma^2}):
N+n-k < \ell < N+n+k-1 \right\}$},
\nonumber\\
&& \alpha^{k-1} H_{N+n}\left( \frac{\alpha}{\sqrt{2 \sigma^2}} \right)
=
\left(\frac{\sigma^2}{2} \right)^{(k-1)/2} H_{N+n+k-1}
\left( \frac{\alpha}{\sqrt{2 \sigma^2}} \right)
\nonumber\\
&& \quad + \mbox{linear combination of
$\left\{H_{\ell}
(\alpha/\sqrt{2 \sigma^2}):
N+n-k < \ell < N+n+k-1 \right\}$},
\nonumber
\end{eqnarray}
for $k=2,3, \dots$.
Therefore, by the multilinearlity of determinant,
$$
\det [A_{N+n-1} \quad A_{N+n} ]=
\sigma^{n(n-1)} \prod_{\ell=2}^{n}
\frac{(N+n-1)!}{(N+n-\ell)!}
\det[B \, C],
$$
where $B$ and $C$ are $2n \times n$ rectangular
matrices whose $(j,k)$-elements are
given by
$$
H_{N+n-k} \left( \frac{\alpha_j}{\sqrt{2 \sigma^2}} \right)
\quad \mbox{and} \quad
H_{N+n+k-1} \left( \frac{\alpha_j}{\sqrt{2 \sigma^2}} \right),
$$
respectively.
Since 
$\det[B \, C]=(-1)^{[n/2]} \det_{1 \leq j, k \leq 2n}
[H_{N+k-1}(\alpha_j/\sqrt{2 \sigma^2})]$,
$(-1)^{[n/2]+n(n-1)/2}=1$ for any $n \in \N$,
and
\begin{equation}
h_n(\alpha_1, \dots, \alpha_n)
h_n(\alpha_{n+1}, \dots, \alpha_{2n})
\prod_{j=1}^{n} \prod_{k=1}^{n}(\alpha_{n+k}-\alpha_j)
=h_{2n}(\alpha_1, \dots, \alpha_{2n}),
\label{eqn:hhh}
\end{equation}
proof of the equivalence of (\ref{eqn:chara1}) and
(\ref{eqn:chara2}) is completed. 
\vskip 0.5cm

Next we prove that 
(\ref{eqn:chara3}) is equal to
\begin{equation}
\frac{\Gamma(N+\nu+n)}{(N+n)!}
\frac{1}{h_{2n}(\valpha)}
\det_{1 \leq j, k \leq 2n}
\left[\widehat{L}^{\nu}_{N+j-1}(\alpha_k; \sigma^2) \right]
\label{eqn:Mnu2}
\end{equation}
for any $\nu > -1$.
By applying the identity (\ref{eqn:Ishikawa1}), we have
\begin{eqnarray}
&& \det_{1 \leq j, k \leq n}
\left[ \frac{1}{\alpha_j-\alpha_{n+k}}
\left| \begin{array}{ll}
L^{\nu}_{N+n}(\alpha_j/2 \sigma^2) &
L^{\nu}_{N+n}(\alpha_{n+k}/2 \sigma^2) \cr
L^{\nu}_{N+n-1}(\alpha_j/2 \sigma^2) &
L^{\nu}_{N+n-1}(\alpha_{n+k}/2 \sigma^2) 
\end{array} \right| \right]
\nonumber\\
&=& 
\frac{(-1)^{n(n-1)/2}}
{\prod_{j=1}^{n} \prod_{k=1}^{n}(\alpha_{n+k}-\alpha_{n})}
\det [\widetilde{A}_{N+n-1} \quad \widetilde{A}_{N+n} ],
\nonumber
\end{eqnarray}
where $\widetilde{A}_{\ell}$ is the $2n \times n$ rectangular
matrix whose $(j,k)$-element is
given by $\alpha_j^{k-1} L^{\nu}_{\ell}
(\alpha_j/2 \sigma^2)$,
$1 \leq j \leq 2n, 1 \leq k \leq n$.
Here we use the recurrence relation of Laguerre
polynomials
\begin{equation}
(n+1) L^{\nu}_{n+1}(x)=
(-x+2n+\nu+1) L^{\nu}_{n}(x)
-(n+\nu) L^{\nu}_{n-1}(x).
\label{eqn:rec2}
\end{equation}
Then we find that
\begin{eqnarray}
&& \alpha^{k-1} L^{\nu}_{N+n-1}
\left( \frac{\alpha}{2 \sigma^2} \right)
=
(-2 \sigma^2)^{k-1} \prod_{\ell=1}^{k-1} (N+n+\nu-\ell) 
L^{\nu}_{N+n-k}
\left( \frac{\alpha}{2 \sigma^2} \right)
\nonumber\\
&& \quad + \mbox{linear combination of
$\Big\{L^{\nu}_{\ell}
(\alpha/2 \sigma^2):
N+n-k < \ell < N+n+k-1 \Big\}$},
\nonumber\\
&& \alpha^{k-1} L^{\nu}_{N+n}
\left( \frac{\alpha}{2 \sigma^2} \right)
=
(-2 \sigma^2)^{k-1} \prod_{\ell=1}^{k-1} (N+n+\ell)
L^{\nu}_{N+n+k-1}
\left( \frac{\alpha}{2 \sigma^2} \right)
\nonumber\\
&& \quad + \mbox{linear combination of
$\Big\{L^{\nu}_{\ell}
(\alpha/2 \sigma^2):
N+n-k < \ell < N+n+k-1 \Big\}$},
\nonumber
\end{eqnarray}
for $k=2,3, \dots$.
Therefore, by the multilinearlity of determinant,
$$
\det [\widetilde{A}_{N+n-1} \quad \widetilde{A}_{N+n} ]=
(2 \sigma^2)^{n(n-1)} 
\left( \frac{\Gamma(N+n+\nu)}{(N+n)!} \right)^{n-1}
\prod_{\ell=1}^{n-1} \frac{(N+n+\ell)!}{\Gamma(N+\nu+\ell)}
\det[\widetilde{B} \, \widetilde{C}],
$$
where $\widetilde{B}$ and $\widetilde{C}$ are $2n \times n$ rectangular
matrices whose $(j,k)$-elements are
given by
$$
L^{\nu}_{N+n-k} \left( \frac{\alpha_j}{2 \sigma^2} \right)
\quad \mbox{and} \quad
L^{\nu}_{N+n+k-1} \left( \frac{\alpha_j}{2 \sigma^2} \right),
$$
respectively.
Since 
$\det[\widetilde{B} \, \widetilde{C}]
=(-1)^{[n/2]} \det_{1 \leq j, k \leq 2n}
[L^{\nu}_{N+k-1}(\alpha_j/2 \sigma^2)]$,
and (\ref{eqn:hhh}), 
the equivalence of (\ref{eqn:chara3}) and
(\ref{eqn:Mnu2}) is proved for $\nu > -1$.
For $\nu \in \N_0$, (\ref{eqn:chara4}) is immediately obtained.
The first expressions in 
(\ref{eqn:chara5}) and (\ref{eqn:chara6}) are obtained
by setting $\nu=1/2$ and $\nu=-1/2$
in (\ref{eqn:Mnu2})
and by applying the relations (\ref{eqn:M3}) and (\ref{eqn:M4}).
For the second expressions 
in (\ref{eqn:chara5}) and (\ref{eqn:chara6}),
we use the following
relations between the Hermite polynomials
and Laguerre polynomials with $\nu=\pm 1/2$,
\begin{eqnarray}
\label{eqn:L1/2}
&& x L_n^{1/2} \left( \frac{x^2}{2} \right)
=\frac{(-1)^n}{ 2^{2n+1/2} n!}
H_{2n+1} \left(\frac{x}{\sqrt{2}} \right),\\
\label{eqn:L-1/2}
&& L_n^{-1/2}\left(\frac{x^2}{2} \right)
=\frac{(-1)^n}{ 2^{2n} n!} 
H_{2n}\left(\frac{x}{\sqrt{2}} \right).
\end{eqnarray} 
Then the proof of Theorem \ref{Theorem:Main2}
is completed. \qed

\vskip 0.5cm
\begin{small}
\noindent{\bf Acknowledgements} \quad
The present author would like to thank N. Minami
for giving him an opportunity to publish this paper
in the RIMS K\^oky\^uroku.


\end{small}
\end{document}